\documentclass[11pt,oneside]{article}
\usepackage{amsfonts}
\usepackage[mathscr]{euscript}
\usepackage{mathrsfs}
\usepackage{xcolor}

\newcommand{\lbl}[1]{\label{#1}}

\newtheorem{theorem}{Theorem}[section]
\newtheorem{proposition}{Proposition}[section]
\newtheorem{lemma}{Lemma}[section]
\newtheorem{remark}{Remark}[section]
\newtheorem{corollary}{Corollary}[section]

\newcommand\bes{\begin{eqnarray}}
\newcommand\ees{\end{eqnarray}}
\newcommand{\bess}{\begin{eqnarray*}}
\newcommand{\eess}{\end{eqnarray*}}
\newcommand{\ds}{\displaystyle}
\newcommand{\dd}{\displaystyle}
\newcommand{\bbb}{\big}
\newcommand\lm{\lambda}
\newcommand\ty{\infty}
\newcommand\vp{\varphi}

\begin{document}
  \pagestyle{myheadings}
\thispagestyle{empty}

\begin{center}{\Large\bf  A reaction-diffusion-advection equation}\\[1mm]
{\Large\bf with mixed and free boundary conditions}\footnote{This work was
supported by NSFC Grant 11371113}\\[4mm]

{\Large Yonggang Zhao}\\[2mm]
Department of Mathematics, Harbin Institute of Technology, Harbin 150001, PR China;\\
College of Mathematics and Information Science, Henan Normal University, Xinxiang
453007, PR China \\[3mm]
{\Large Mingxin Wang}\footnote{Corresponding author: mxwang@hit.edu.cn}\\[1mm]
{Natural Science Research Center, Harbin Institute of Technology, Harbin 150080, PR China}
\end{center}

\begin{quote}
\noindent{\bf Abstract.} We investigate a reaction-diffusion-advection equation of the form $u_t-u_{xx}+\beta u_x=f(u)$ $(t>0,\,0<x<h(t))$ with mixed boundary condition at $x=0$ and a free boundary condition at $x=h(t)$. Such a model may be applied to describe the dynamical process of a new or invasive species adopting a combination of random movement and advection upward or downward along the resource gradient, with the free boundary representing the expanding front. The goal of this paper is to understand the effect of advection environment and no flux across the left boundary on the dynamics of this species. When $|\beta|<c_0$, we first derive the spreading-vanishing dichotomy and sharp threshold for spreading and vanishing. Then provide a much sharper estimate for the spreading speed of $h(t)$ and the uniform convergence of $u(t,x)$ when spreading happens. For the case $|\beta|\geq c_0$, some results concerning spreading, virtual spreading, vanishing and virtual vanishing are obtained. Where $c_0$ is the minimal speed of traveling waves of the differential equation.

\noindent{\bf Keywords:} Reaction-diffusion-advection equation; free boundary; spreading and vanishing; sharp threshold; long time behavior.

\noindent {\bf AMS subject classifications (2010)}:
35K20, 35R35, 92B05, 35B40.
 \end{quote}

 \section{Introduction}
 \setcounter{equation}{0} {\setlength\arraycolsep{2pt}

In recent years there has been growing interest in understanding the role that the free boundary plays in the dynamics of introduction of beneficial species or invasion of harmful species. In reality the dispersal of new or invasive species is often nonrandom as both dispersal rate and direction can depend upon a combination of local biotic and abiotic factors such as climate, food, and conspecifics. For instance, some diseases spread along the wind direction. From a mathematical point of view, to take into account the influence of advection,  one of the simplest but probably still realistic approaches  is to assume that species can move upward or downward along the gradient of the density (see, for example, \cite{BH,HL,RSB,SG,VL}).

Recently, Gu, Lin \& Lou \cite{GLL,GLL2}, Kaneko \& Matsuzawa \cite{KM}, and Gu, Lou \&  Zhou \cite{GLZ} studied the influence of positive advection on the long time behavior of solutions to the following problem:
\begin{eqnarray}\label{1.1}
\left\{\begin{array}{ll}
u_t-u_{xx}+\beta u_x=f(u), &t>0,\;g(t)<x<h(t),\\[1mm]
u(t,g(t))=0,\ \; g'(t)=-\mu u_x(t,g(t)),&t\geq0,\\[1mm]
u(t,h(t))=0,\ \; h'(t)=-\mu u_x(t,h(t)), \ \ &t\geq0,\\[1mm]
-g(0)=h_0=h(0),\ \; u(0,x)=u_0(x),& -h_0\leq x\leq h_0,
\end{array}
\right.
\end{eqnarray}
where $\beta$, $\mu$ and $h_0$ are positive constants, $u_0$ is a nonnegative $C^2$ function with the support on $[-h_0,h_0]$. Problem (\ref{1.1}) arises in modeling the spreading of a new or invasive species going through the influence of dispersal and advection (expressed by $\beta u_x$). The unknown $u(t,x)$ represents the population density over a one dimensional habitat and the free boundaries $x=g(t)$ and $x=h(t)$ stand for the expanding fronts of the species.

In \cite{GLL,GLL2}, the authors considered the asymptotic behavior of solutions to (\ref{1.1}) when the advection coefficient $\beta\in(0,2)$ and $f(u)=u(1-u)$. They obtained a spreading-vanishing dichotomy, namely the solution either converges to 1 locally uniformly in $\mathbb{R}$ or to $0$ uniformly in its occupying domain. Moreover, by introducing a parameter $\sigma$ in the initial value, they exhibited a sharp threshold between spreading and vanishing, that is, there exists a nonnegative constant $\sigma^*$ such that spreading happens if $\sigma>\sigma^*$, and vanishing happens if $\sigma\leq\sigma^*$.  Furthermore, they derived the following conclusions concerning the asymptotic spreading speed: if spreading happens, then there exist two positive constants $c^*_l$ and $c^*_r$ with $c^*_l<c^*_r$ such that
   $$\lim_{t\to\infty}\frac{g(t)}t=-c^*_l,\ \ \lim_{t\to\infty}\frac{h(t)}t=c^*_r.$$
For the general case that $f(u)$ is monostable, bistable or of combustion type, the above result is improved in \cite{KM} to a much sharper estimate for the different spreading speeds of the fronts: when $\beta\in(0,c_0)$,
$$\lim_{t\to\infty}g'(t)=-c^*_l,\ \ \lim_{t\to\infty}h'(t)=c^*_r,\ \ \lim_{t\to\infty}[g(t)+c^*_lt]=G_\infty,\ \ \lim_{t\to\infty}[h(t)-c^*_rt]=H_\infty$$
for some $G_\infty,H_\infty\in\mathbb{R}$, where $c_0$ is the minimal speed of the traveling waves:
 \begin{eqnarray}\label{1.2}
\left\{\begin{array}{ll}
q''-cq'+f(q)=0, \ \ q>0 \ \ \ {\rm in}\ \ \mathbb R,\\[1mm]
q(-\ty)=0, \ \ q(\ty)=1.
\end{array}
\right.
  \end{eqnarray}
The reader interested in the number $c_0$ can refer to \cite{AW2,AW,DLou}. Apart from the above result, the authors of \cite{KM} described how the solution approaches a semi-wave when the nonlinear function is a monostable, bistable or combustion type. Very recently, for the general function $f(u)\in C^1([0,\ty))$  satisfying the following condition $(F)$,  Gu, Lou \& Zhou \cite{GLZ} extended the advection coefficient $\beta$ to $\beta\in(0,\infty)$. They found a critical value $\beta^*$ with $\beta^*>c_0$ and gained the trichotomy results for $c_0\leq \beta<\beta^*$, vanishing result for $\beta\geq\beta^*$, where $c_0=2\sqrt{f'(0)}$ is the minimal speed of the traveling waves of (\ref{1.2}).

Motivated by the above works, in this paper we concern with the following reaction-diffusion-advection model with a free boundary:
  $$\left\{\begin{array}{ll}
u_t-u_{xx}+\beta u_x=f(u), &t>0,\;0<x<h(t),\\[1mm]
B[u](t,0)=0,\ u(t,h(t))=0,&t\geq0,\\[1mm]
h'(t)=-\mu u_x(t,h(t)),&t\geq0,\\[1mm]
h(0)=h_0,\ \ u(0,x)=u_0(x),& 0\leq x\leq h_0,
  \end{array} \right.\eqno{(P)}$$
where $\beta\in\mathbb{R}$, and $\mu,\,h_0>0$; the left boundary operator $B[u]=au-bu_x$ with $a,b\geq0$ and $a+b>0$;  $x=h(t)$ represents the moving boundary which is to be determined with the solution $u(t,x)$; $f: [0,\infty)\to\mathbb{R}$ is a $C^1$  function and satisfies
 $$\left\{\begin{array}{ll}
f(0)=0=f(1), \ \ (1-u)f(u)>0\ \ \ {\rm for\ } u>0,\ u\neq1,\\[1mm]
f'(0)>0,\ \ f'(1)<0,\ \ f(u)\leq f'(0)u\ \ \ {\rm for\ } u>0;
\end{array}\right.\eqno{(F)}$$
the initial function $u_0$ belongs to $\mathscr{X}(h_0)$, where
\begin{eqnarray*}
\mathscr{X}(h_0)=\left\{\psi\in C^2([0,h_0]):B[\psi](0)=\psi(h_0)=0,\ \psi'(h_0)<0,\ \psi(x)>0 \ {\rm in}\ (0,h_0)
\right\}.
\end{eqnarray*}

For more clarity, in the following we always denote
 $$c_0=2\sqrt{f'(0)},$$
which is the minimal speed of the traveling waves of (\ref{1.2}).

The main intention of this paper is to study the dynamics of the problem $(P)$ under the assumptions that the advection coefficient is a real number not only a positive one and there is a mixed boundary condition at the left boundary $x=0$. When $|\beta|<c_0$, we will provide a rather complete description of the spreading-vanishing dichotomy, sharp threshold for spreading and vanishing, sharp asymptotic spreading speed and the uniform convergence of the solution when spreading happens.  Moreover, we will briefly describe the long time behavior of solutions in cases either $\beta\geq c_0$ and $a\geq bc_0/2$ or $\beta\leq- c_0$.

Applying an analogous argument as in \cite[Theorem 2.1]{Wang15, Wcnsns15} one can demonstrate that, for any given $u_0\in\mathscr{X}(h_0)$, problem $(P)$ admits a unique time-global solution $(u,h)$ with $u\in C^\infty((0,\infty)\times[0,h(t)])$ and $h\in C^\infty((0,\infty))$. Moreover, for any $\alpha\in(0,1)$, there exists a positive constant $C$ depending on $\beta, \alpha, h_0, a,b$ and $\|u_0\|_\infty$, such that $0<u(t,x)\leq C$, $0<h'(t)\leq C$ for all $ t>0$ and $0<x<h(t);$ and
\begin{eqnarray}
\|u(t,\cdot)\|_{C^1([0,h(t)])}\leq C,\ \forall \ t\geq1;\ \ \|h'\|_{C^{\frac{1+\alpha}{2}}([n+1,n+3])}\leq C, \ \forall\ n\geq0. \label{1.3}
\end{eqnarray}
Denote $h_\infty=\lim_{t\to\infty}h(t)$ for simplicity.

The stationary problem for $(P)$ can be written as
\begin{eqnarray}\label{1.4}
\left\{\begin{array}{l}
v''-\beta v'+f(v)=0, \ \ 0<x<\ell,\\[1mm]
B[v](0)=0
\end{array}
\right.
\end{eqnarray}
for some $0<\ell\leq\infty$. By the phase plane analysis (see \cite{AW,DLou,GLZ}), it is not difficult to derive that the non-negative solutions of problem (\ref{1.4}) fall into the following categories:

(i) constant solution: $v=0$ or $v=1$ (the latter can exist when $a=0$);

(ii)  strictly increasing solution on the half-line: $v(x)=\tilde v(x)$, where $\tilde v(x)$ satisfies (\ref{1.4}) on $[0,\infty)$, $\tilde v'(x)>0$ and $\tilde v(\infty)=1$ (such a solution exists uniquely when $\beta<c_0:=2\sqrt{f'(0)}$ and $a>0$);

(iii) non-monotone solution on the half-line: $v(x)=\check{v}(x)$, where $\check{v}(x)$ satisfies (\ref{1.4}) on $[0,\infty)$, and $\check{v}(\infty)=0$ (such a solution exists uniquely when $\beta\leq -c_0$);

(iv) solutions with finite interval length $\ell$: $v(x)=v_\ell(x)$, where $v_\ell(x)$ satisfies (\ref{1.4}) on $[0,\ell)$, and $v_\ell(\ell)=0$ (such solutions can exist when $-c_0<\beta<c_0$).

When $\beta=0$ (without the advection term), problem $(P)$ becomes
 \bes\left\{\begin{array}{ll}
u_t-u_{xx}=f(u), &t>0,\;0<x<h(t),\\[1mm]
B[u](t,0)=0,\ u(t,h(t))=0,&t\geq0,\\[1mm]
h'(t)=-\mu u_x(t,h(t)),&t\geq0,\\[1mm]
h(0)=h_0,\ \ u(0,x)=u_0(x),& 0\leq x\leq h_0.
  \end{array} \right.\lbl{a.1}\ees
This problem and other related problems have been studied by many authors. Recently, Liu \& Lou \cite{LL, LL1} studied problem (\ref{a.1}) when $a>0$ and $f(u)$ is a monostable, bistable or combustion type nonlinearity. Du \& Lin \cite{DL} investigated problem (\ref{a.1})  initially with $a=0$ ($B[u](t,0)=u_x(t,0)=0$) for the logistic nonlinearity $f(u)=u(1-u)$; while Kaneko \& Yamada \cite{KY} studied problem (\ref{a.1}) with $b=0$ and $f(u)=u(1-u)$ or $f(u)=u(u-c)(1-u)$. When the nonlinear term $f(u)$ is replaced by $u(\alpha(x)-u)$ and $\alpha(x)$ changes signs, problem (\ref{a.1}) was studied by Wang
\cite{Wang15}. When the nonlinear term $f(u)$ is replaced by $u(\alpha(t,x)-\beta(t,x)u)$ and $\alpha(t,x)$, $\beta(t,x)$ are $T$-periodic in time $t$, problem (\ref{a.1}) was studied by Du, Guo \& Peng \cite{DGP} (with $a=0$) and Wang \cite{Wperiodic15} ($\alpha(t,x)$ changes signs).
Du \& Guo \cite{DG, DG1} and Du \& Liang \cite{DLiang} considered the higher dimension and heterogeneous environment case.  Peng \& Zhao \cite{PZh} studied one seasonal succession case. When the nonlinear term $f(u)$ is replaced by a general function including monostable, bistable and combustion type, the double free boundary problems has been considered by Du \& Lou \cite{DLou}, Du, Matsuzawa \& Zhou \cite{DMZ} and Kaneko \cite{Ka}. The diffusive competition system with a free boundary has been researched by Du \& Lin \cite{DL2}, Guo \& Wu \cite{GW, GW15}, Wang \& Zhao \cite{WZh1} and Wu \cite{Wu13, Wu14}. The diffusive prey-predator model with free boundaries has been studied by Wang \&  Zhao \cite{WMX, Wcnsns15, ZW}. Recently, Ge et al. \cite{GKLZ} and Huang \& Wang \cite{HW} investigated the epidemic model with a free boundary, Li \& Lin  \cite{LiLin} discussed a mutualistic model with advection and a free boundary.

The organization of this paper is as follows. Sections \ref{se2} and \ref{se3} deal with the case $|\beta|<c_0$. Section \ref{se4} concerns the case either $\beta\geq c_0$ and $a\geq bc_0/2$, or $\beta\leq- c_0$

In Section \ref{se2}, the spreading-vanishing dichotomy and two kinds of sharp thresholds between spreading and vanishing are displayed for $|\beta|<c_0$. To do this, we introduce two eigenvalue problems, and then find the critical interval width $\ell^*$ of free boundary $h(t)$, that is, $h_\infty<\infty$ implies $h_\infty\leq \ell^*$, and $h_0\ge \ell^*$ implies $h_\infty=\infty$. Moreover, in order to establish the sharp thresholds for spreading and vanishing more perfectly, we construct three groups of upper solutions to show that the sharp thresholds are positive, see Theorems \ref{th2.5} and \ref{th2.6}.

In Section \ref{se3}, we concern the long time behavior of $(h,u)$ for the spreading case and $|\beta|<c_0$. A much sharper estimate for the spreading speed o $h(t)$, and the uniform convergence of $u(t,x)$ will be given. The approach for the sharper estimate of spreading speed is to use a zero number argument inspired by \cite{KM}, utilizing a one-dimensional problem for a single equation. Moreover, the uniform convergence of $u(t,x)$ is obtained by combining two groups of upper and lower solutions with the locally uniform convergence near two boundaries $x=0$ and $x=h(t)$. Though the outline of the approach in this part largely follows that of \cite{GLL, GLL2, KM, GLZ}, some of the technical proofs here are different from and much more involved than the corresponding ones in those references. Our argument involves some new ideas and techniques.

In Section \ref{se4} we deal with the case either $\beta\geq c_0$ and $a\geq bc_0/2$, or $\beta\leq- c_0$. Firstly, it is proved that $u(t,\cdot)$ converges to $0$ locally uniformly in $[0,h_\infty)$ as $t\to\ty$ regardless of $h_\infty<\infty$ or $h_\infty=\infty$. When either $\beta\geq\beta^*$ and $a\geq bc_0/2$, or $\beta<-c_0$, we how that $h_\infty<\infty$. Moreover, some results concerning spreading, virtual spreading, vanishing and virtual vanishing are obtained.

We should remark that, to our knowledge, the present paper is the first one  investigating the negative advection case ($\beta<0$).

\section{Spreading-vanishing dichotomy and sharp threshold}\label{se2}
\setcounter{equation}{0}

In this section we first provide two types of comparison principles and some definitions of upper solutions, then give the spreading-vanishing dichotomy.  In the third subsection we establish two sharp thresholds for spreading and vanishing, which depend on different varying parameters. One is $\mu$, and the other is the initial value $u_0(x)$.

\subsection{Two comparison principles}
 In this subsection we provide two types of comparison principles which are an analogue of \cite[Lemmas 3.5]{DL} and the proofs will be omitted.

\begin{lemma}\label{l2.1}{\hspace{-1.8mm}\bf .}
Let $(u,h)$ be a solution of $(P)$. Assume that $(\bar u,\bar h)\in C(\overline \mathcal{D})\cap C^{1,2}(\mathcal{D})\times C^1([0,\infty))$ with $\mathcal{D}=\{(t,x):\,t>0, 0< x<\bar h(t)\}$, satisfies
\begin{eqnarray*}
\left\{\begin{array}{lll}
 \bar u_t- \bar u_{xx}+\beta\bar u_x\geq f(\bar u), &t>0, \   0<x<\bar h(t),\\[1mm]
B[\bar u](t,0)\geq0,\ \bar u(t,\bar h(t))=0,\ \ &t\geq0,\\[1mm]
\bar h'(t)\geq-\bar\mu \bar u_x(t,\bar h(t)),\ \ &t\geq0.
 \end{array}\right.
\end{eqnarray*}
If $\bar\mu\ge\mu$, $\bar h(0)\geq h_0$, $\bar u(0,x)\geq u_0(x)$  for all $0\leq x\leq h_0$, then
$$h(t)\leq\bar h(t),\ \ \forall\ t\geq0;\quad\ u(t,x)\leq\bar u(t,x),\ \  \forall\ t\geq0,\ 0\leq x\leq h(t).$$
\end{lemma}

\begin{lemma}\label{l2.2}{\hspace{-1.8mm}\bf .}
Let $(u,h)$ be a solution of $(P)$. Assume that $\bar g,\bar h\in C^1([0,\infty))$ and $0\le\bar g(t)<\bar h(t)$, $\bar g(t)\leq h(t)$ for all $t\ge 0$, $\bar u\in C(\overline O)\cap C^{1,2}(O)$ with $O=\{(t,x):\,t>0,\,\bar g(t)< x<\bar h(t)\}$. If
\begin{eqnarray*}
\left\{\begin{array}{lll}
 \bar u_t- \bar u_{xx}+\beta\bar u_x\geq f(\bar u), &t>0, \  \bar g(t)<x<\bar h(t),\\[1mm]
\bar u(t,\bar g(t))\geq u(t,\bar g(t)),\ \ &t\geq0,\\[1mm]
\bar u(t,\bar h(t))=0,\ \bar h'(t)\geq-\bar\mu u_x(t,\bar h(t)),\ \ &t\geq0,
 \end{array}\right.
\end{eqnarray*}
and $\mu\le\bar\mu$, $h_0\leq\bar h(0)$, $u_0(x)\leq \bar u(0,x)$ for all $\bar g(0)\leq x\leq h_0$, then
$$h(t)\leq\bar h(t),\ \ \forall\ t\geq0;\quad\ u(t,x)\leq\bar u(t,x)\ \ \forall\ t\geq0,\ \bar g(t)\leq x\leq h(t).$$
\end{lemma}

The pair $(\bar u,\bar h)$ in Lemma \ref{l2.1} and the triple $(\bar u, \bar g,\bar h)$ in Lemma  \ref{l2.2} are usually called an upper solution of problem $(P)$. We can define a lower solution  by reversing all the inequalities in the obvious places. In addition, the corresponding results for lower solutions can be also shown by the similar manner.

\subsection{Spreading-vanishing dichotomy}

Firstly, the asymptotic behavior of $u$ is presented for vanishing case $h_\ty<\ty$. To do this, we give the following more general lemma, whose proof is essentially similar to that of \cite[Theorem 2.2]{WZh1}. We will leave out the details because the advection term, mixed boundary condition and more general reaction term do not influence the availability of the argument in \cite[Theorem 2.2]{WZh1}.

\begin{lemma}\lbl{l2.3}{\hspace{-1.8mm}\bf .} Let $\beta, c\in\mathbb{R}$ and $\mu>0$. Assume that $s\in C^1([0,\infty))$, $w\in C^{\frac{1+\nu}2,1+\nu}([0,\infty)\times[0,s(t)])$ and satisfy $s(t)>0$, $w(t,x)>0$ for $t\ge 0$ and $0<x<s(t)$. We further suppose that
$\ds\lim_{t\to\infty} s(t)<\infty$, $\ds\lim_{t\to\infty} s'(t)=0$ and there exists a constant $C>0$ such that $\|w(t,\cdot)\|_{C^1[0,\,s(t)]}\leq C$ for $t>1$. If $(w,s)$ satisfies
  \bess\left\{\begin{array}{lll}
 w_t-w_{xx}+\beta w_x\geq cw, &t>0, \ 0<x<s(t),\\[.5mm]
 B[w]x=0 ,\ \ \ &t>0, \ \ x=0,\\[.5mm]
 w=0, \ s'(t)\geq-\mu w_x, \ &t\ge 0, \ x=s(t),
 \end{array}\right.\eess
then $\ds\lim_{t\to\infty}\max_{0\leq x\leq s(t)}w(t,x)=0$.
 \end{lemma}

By the second estimate of (\ref{1.3}), we see that if $h_\ty<\ty$, then $\ds\lim_{t\to\ty}h'(t)=0$. As a consequence of Lemma \ref{l2.3} we derive the asymptotic behavior of $u$ for $h_\ty<\ty$.

\begin{theorem}\lbl{th2.1}{\hspace{-1.8mm}\bf .} Let $(u,h)$ be the solution of $(P)$. If $h_\ty<\ty$, then $\dd\lim_{t\to\infty}\max_{0\leq x\leq h(t)}u(t,x)=0$.
 \end{theorem}

Then we provide a locally uniform convergence theorem   for  spreading case $h_\infty=\infty$.

\begin{theorem}\label{th2.2}\hspace{-1mm}{\rm(}Local convergence{\rm)}.
Assume that $|\beta|<c_0$, $(u,h)$ is the solution of problem $(P)$ and $h_\ty=\ty$. Then we have

{\rm (i)} in case $a>0$, $u$ converges to $\tilde v$ locally uniformly in $[0,\infty)$, where $\tilde v(x)$ is the strictly increasing solution of $(\ref{1.4})$ on the half-line;

{\rm(ii)} in case $a=0$,  $u$ converges to $1$ locally uniformly in $[0,\infty)$.
\end{theorem}

\noindent{\bf Proof.}
Applying an analogous argument as in \cite[Theorem 1.1]{DM}, \cite[Theorem 1.1]{DLou} and \cite[Theorem 2.1]{GLL}, one can demonstrate that when $t$ goes to $\infty$, $u(t,x)$ must approach a stationary solution of problem $(P)$, that is, a solution $v$ of (\ref{1.4}), locally uniformly in $[0,\infty)$. Moreover, the sole possible choice for the $\omega$-limit of $u$ in the topology of $L^\infty_{\rm loc}([0,\infty))$ is $0$ or $\tilde v$ when $a>0$, and is $0$ or $1$ when $a=0$.

If the element $0$ does not belong to the $\omega$-limit of $u$ for spreading case, then the desired results of this theorem are deduced immediately.

In the following we shall prove that the $\omega$-limit of $u$ does not include $0$  for  spreading case.
Note that $h_\ty=\ty$, there exists $\tau>0$ such that $h(\tau)>2\pi/\sqrt{c^2_0-\beta^2}$. Let $\ell\in(2\pi/\sqrt{c^2_0-\beta^2},h(\tau))$. Then there exist small positive constants $\sigma$ and $\varepsilon$ such that
$$\frac{4\pi^2}{\ell^2}\leq 4(f'(0)-\sigma)-\beta^2,\ \ f(s)\geq (f'(0)-\sigma)s\ \ {\rm for}\ s\in [0,\varepsilon].$$
Define
  $$\tilde  w(x)=\varepsilon e^{\frac\beta2x}\sin\frac{\pi x}{\ell}, \ \ 0\leq x\leq \ell.$$
It is obvious that $\tilde w(0)=0=\tilde w(\ell)$. An elementary calculation yields, for $0< x<\ell$,
$$-\tilde w''+\beta\tilde w'-f(\tilde w)\leq\tilde w\left(\frac{\beta^2}4+\frac{\pi^2}{\ell^2}-f'(0)+\sigma\right)\leq0.$$
Since $u(t,x)>0$ for all $t>0$ and $0<x<h(t)$, we have $u(t,\ell)>0$ for $t\geq \tau$, and $u(\tau,x)\geq \tilde w(x)$ on $[0,\ell]$ provided $\varepsilon>0$ is sufficiently small. Besides, the boundary condition $B[u](t,0)=0$ implies $u(t,0)\geq0$ for any $t>0$. It follows from the comparison principle that $u(t,x)\geq\tilde w(x)$ for all $t>\tau$ and $0\leq x\leq \ell$. This indicates that $0$ is not in the $\omega$-limit of $u$. The proof is completed.\ \ \ \fbox{}

\vspace{.4cm}
Next we introduce the following two eigenvalue problems
 \bes\left\{\begin{array}{ll}
 -\vp''+\beta\vp'-f'(0)\vp=\zeta\vp, \ \ 0<x<\ell,\\[1mm]
 B[\vp](0)=0, \ \ \vp(\ell)=0
 \end{array}\right.\lbl{2.1}\ees
and
 \bes\left\{\begin{array}{ll}
 -\phi''-f'(0)\phi=\gamma\phi, \ \ 0<x<\ell,\\[1mm]
 B[\phi](0)=0, \ \ \phi(\ell)=0,
 \end{array}\right.\lbl{2.2}\ees
where $\ell>0$ is a constant. Let $\zeta_1(\ell)$ and $\gamma_1(\ell)$  be the first eigenvalues of (\ref{2.1}) and (\ref{2.2}), respectively. By a careful calculation we achieve the following conclusions:

(i)  Both $\zeta_1(\ell)$ and $\gamma_1(\ell)$ are continuous and strictly decreasing in $\ell$;

(ii) $\lim_{\ell\to 0^+}\zeta_1(\ell)=\lim_{\ell\to 0^+}\gamma_1(\ell)=\infty$, $\lim_{\ell\to\ty}\zeta_1(\ell)=\beta^2/4-f'(0)$, $\lim_{\ell\to\ty}\gamma_1(\ell)=-f'(0)$.

By virtue of the above conclusions (i) and (ii), it is easy to see that if $|\beta|<c_0$, i.e., $\beta^2<4f'(0)$, then there exist $\ell^*,\,\ell_*>0$ such that
\begin{equation}\label{zz}\zeta_1(\ell^*)=0, \ \ \gamma_1(\ell_*)=-\beta^2/4. \end{equation}
By the monotonicity of $\zeta_1(\ell)$ and $\gamma_1(\ell)$, we have

(i) \, $\zeta_1(\ell)>0$ for $\ell<\ell^*$, and $\zeta_1(\ell)<0$ for $\ell>\ell^*$;

(ii)\, $\gamma_1(\ell)>-\beta^2/4$ if $\ell<\ell_*$, while $\gamma_1(\ell)<-\beta^2/4$ when  $\ell>\ell_*$.\\
For the case $b=0$, we can compute that $\zeta_1(\ell)=\beta^2/4+\pi^2/\ell^2-f'(0)$ and $\gamma_1(\ell)=\pi^2/\ell^2-f'(0)$, and that $\ell^*$, $\ell_*$ satisfying (\ref{zz}) are both equal to $2\pi/\sqrt{c_0^2-\beta^2}$.

\vskip 4pt
The following lemma is an analogue of \cite[Lemma 3.2]{Wang15}, so the details are omitted here.

 \begin{lemma}\lbl{l2.4}{\hspace{-1.8mm}\bf .} If $h_\infty<\infty$, then $\zeta_1(h_\infty)\geq 0$.
\end{lemma}

Taking advantage of Lemma \ref{l2.4} we can achieve

\begin{theorem}\lbl{th2.3}{\hspace{-1.8mm}\bf .} Let $|\beta|<c_0$ and $(u,h)$ be the solution of $(P)$. If $h_\ty<\ty$, then $h_\ty\le\ell^*$. Therefore, $h_0\ge\ell^*$ implies $h_\ty=\ty$ for any $u_0\in \mathscr{X}(h_0)$ and $\mu>0$.
 \end{theorem}

Combining Theorems \ref{th2.1}, \ref{th2.2} and \ref{th2.3}, we derive the main conclusion of this subsection:

\begin{theorem}\label{th2.4}\hspace{-1mm}{\rm(}Spreading-vanishing dichotomy{\rm)}.
Assume that $|\beta|<c_0$ and $(u,h)$ is a solution of $(P)$. Then either

{\rm(i)} vanishing: $h_\infty\le\ell^*$  and $\lim_{t\to\infty}\|u(t,\cdot)\|_{C([0,h(t)])}=0$;
or

{\rm(ii)} spreading: $h_\infty=\infty$, and in case $a>0$, $\lim_{t\to\infty}u(t,x)=\tilde v(x)$ locally uniformly in $[0,\infty)$; in case $a=0$, $\lim_{t\to\infty}u(t,x)=1$ locally uniformly in $[0,\infty)$. Here $\tilde v(x)$ is is the strictly increasing solution of $(\ref{1.4})$ on the half-line.
\end{theorem}
\subsection{Sharp threshold for spreading and vanishing}
In this subsection we first discuss $\mu$ as varying parameters to describe the threshold for spreading and vanishing, then establish the threshold on the initial value, which separates vanishing and spreading.


\begin{theorem}\hspace{-1mm}{\rm(}Threshold on $\mu${\rm)}.\label{th2.5}
Let $|\beta|<c_0$ and $(u,h)$ be any solution of $(P)$. Assume that one of the following conditions holds:

${\rm(i)}$\ $h_0<\ell^*$, and either $b=0$, or $b>0$ and $\beta\le 0$;

  ${\rm(ii)}$\ $h_0<\pi/\sqrt{c^2_0-\beta^2}$, and either $\beta\leq0$, or $\beta>0$ and $a\geq b\beta/2$.\\
Then for any given $u_0\in \mathscr{X}(h_0)$, there exists $\mu^*>0$
such that $h_\infty=\infty$ for $\mu>\mu^*$, while $h_\infty<\ty$ for $\mu\leq\mu^*$.
\end{theorem}

The proof of Theorem \ref{th2.5} will be divided into the following three lemmas: Lemmas \ref{lb.5}-\ref{lb.7}. For the later use, we first give a proposition

\begin{proposition}\lbl{p2.1}{\hspace{-1.8mm}\bf .} Let $C>0$ be a constant. For any given constants $\bar h_0, H>0$, and any given function $\bar u_0\in \mathscr{X}(h_0)$, there exists $\mu^0>0$, depending on $\beta$, $C$, $\bar u_0(x)$ and $\bar h_0$, such that when $\mu\geq\mu^0$ and $(\bar u, \bar h)$ satisfies
 \bess
 \left\{\begin{array}{ll}
   \bar u_t-\bar u_{xx}+\beta\bar u_x\geq -C \bar u, \ &t>0, \ 0<x< \bar h(t),\\[1mm]
  B[\bar u](t,0)=0=\bar u(t, \bar h(t)),\ &t\geq 0,\\[1mm]
 \bar h'(t)=-\mu \bar u_x(t, \bar h(t)), \ &t\geq 0,\\[1mm]
  \bar h(0)=\bar h_0, \ \bar u(0,x)=\bar u_0(x),\ &0\leq x\leq \bar h_0,
 \end{array}\right.
 \eess
we must have $\lim_{t\to\infty}\bar h(t)>H$.
\end{proposition}

This proposition can be proved by the similar method to that of \cite[Lemma 3.2]{WZh1} and the details will be omitted  here.

\begin{lemma}\lbl{lb.5}{\hspace{-1.8mm}\bf .} Let $|\beta|<c_0$ and $(u,h)$ be the solution of $(P)$. If $h_0<\ell^*$ and $b=0$, then the conclusions of Theorem $\ref{th2.5}$ hold.
 \end{lemma}

 \noindent{\bf Proof.} Notice $\ell^*=\ell_*=2\pi/\sqrt{c_0^2-\beta^2}$ for the case $b=0$. The proof consists of three steps.

{\it Step 1}. We prove that for any given $u_0\in \mathscr{X}(h_0)$, there exists $\mu_0>0$ depending on $\beta,h_0,f'(0)$ and $u_0(x)$ such that $h_\infty<\ty$ for $\mu\leq\mu_0$.

One can know $\gamma_1:=\gamma_1(h_0)>-\beta^2/4$ if $h_0<\ell_*$.  Let $\phi$ be the positive eigenfunction of (\ref{2.2}) corresponding to $\gamma_1$. Noting that $\phi'(h_0)<0$, $\phi(0)>0$ when $b>0$, and $\phi'(0)>0$ when $b=0$, it is easy to see that there exists $k>0$ such that
 \bes
 x\phi'(x)\le k\phi(x) , \ \ \forall \  \ 0\le x\le h_0.
 \lbl{2.3}\ees
Let $0<\delta<1$ and $K>0$ be constants, which will be determined later. Set
 \bess
 \dd g(t)=1+\delta-\frac\delta 2{\rm e}^{-\delta t}, \ \
 v(t,x)=Ke^{-\delta t}e^{-\frac\beta2(h_0g(t)-x)}\phi\left(x/g(t)\right), \ \ \ t\geq 0,\ \ 0\leq x\leq h_0g(t).
 \eess
Denote $y=x/g(t)$. Owing to the inequality (\ref{2.3}), $\gamma_1>-\beta^2/4$ and $f(v)\leq f'(0)v$, by routine calculations we show
 \bes
 v_t-v_{xx}+\beta v_x-f(v)&=&v\left(-\delta-\frac{\beta}{2}h_0g'(t)-\frac{yg'(t)\phi'(y)}{g(t)\phi(y)}+\frac{\beta^2}4-\frac{\phi''(y)}
 {\phi(y)g^2(t)}-\frac{f(v)}v\right) \nonumber\\[.5mm]
 &\geq&v\left(-\delta-\frac{\beta}4h_0\delta^2e^{-\delta t}-\frac{y\phi'(y)\delta^2}{\phi(y)g(t)}e^{-\delta t}+\frac{\beta^2}4+\frac{f'(0)+\gamma_1}{g^2(t)}-f'(0)\right)\nonumber\\[.5mm]
 &\geq& v\left(-\delta-\frac{\beta}{4}h_0\delta^2-k\delta^2+\frac{\beta^2}4
 +\frac{\gamma_1}{(1+\delta)^2}-\frac{\delta(2+\delta)}{(1+\delta)^2}f'(0)\right)\nonumber\\
 &>&0, \ \ \forall \ t>0, \ \ 0<x<g(t)\qquad
 \lbl{2.4}\ees
provided $0<\delta\ll 1$. It is easy to see that
\begin{equation}\label{2.5}
B[v](t,0)=0,\ \ \ v(t,h_0g(t))=0,\ \ \forall\  t\geq0
\end{equation}
as $b=0$. Fix $0<\delta\ll1$. By the regularities of $u_0(x)$ and $\phi(x)$ we can select $K\gg 1$ such that
\begin{equation}\label{2.6}
u_0(x)\leq Ke^{-\frac\beta2(h_0(1+\frac\delta2)-x)}\phi\left( 2x/{(2+\delta)}\right)=v(0,x), \ \ \forall\ 0\leq x\leq h_0.
\end{equation}
Due to $v_x(t,h_0g(t))=K{\rm e}^{-\delta t}\phi'(h_0)/g(t)$ and $\phi'(h_0)<0$, there exists $\mu_0>0$ such that
 \begin{equation}\label{2.7}
h_0g'(t)=h_0\delta^2 {\rm e}^{-\delta t}\geq -\mu v_x(t,h_0g(t)), \ \ \forall\ 0<\mu\leq\mu_0, \ t\geq0.
 \end{equation}

Combining (\ref{2.4})-(\ref{2.7}) we can make use of the comparison principle (Lemma \ref{l2.1}) to achieve that
$h(t)\leq h_0g(t)$ for all $t\geq0$. Thus $h_\infty\leq h_0(1+2\delta)$ for all $0<\mu\leq \mu_0$.

\vskip 4pt
{\it Step 2}. Let $C=\ds\max_{0\le z\le\|u\|_\ty}|f'(z)|$. Then we have $f(u)\ge -Cu$. Take $H=\ell^*$, then it follows from Proposition \ref{p2.1} that there exists $\mu^0>0$ such that $h_\infty>\ell^*$ for $\mu>\mu^0$.
This implies $h_\infty=\ty$ for $\mu>\mu^0$ by Theorem \ref{th2.3}.

{\it Step 3}. Based upon the results obtained by Steps 1 and 2, by use of the continuity method we can prove the the desired result. Please refer to the proof of \cite[Theorem 3.9]{DL} for details.
\ \ \ \fbox{}

\begin{lemma}\lbl{lb.6}{\hspace{-1.8mm}\bf .} Let $|\beta|<c_0$ and $(u,h)$ be the solution of $(P)$. If $h_0<\ell^*$, and $b>0$, $\beta\le 0$, then the conclusions of Theorem $\ref{th2.5}$ hold.
 \end{lemma}

\noindent{\bf Proof.} The proof is essentially similar to that of Lemma \ref{lb.5}, but the difference is the selection of auxiliary function. Hence we only manifest the following
  \vspace{-1mm}\begin{quote}{\bf Claim}: for any given $u_0\in \mathscr{X}(h_0)$, there exists $\mu_0>0$ depending on $\beta,h_0,\zeta_1(h_0)$ and $u_0(x)$ such that $h_\infty<\ty$ for $\mu\leq\mu_0$,\vspace{-1mm}\end{quote}
and the rest proof will be omitted.

Since $h_0<\ell^*$, we have $\zeta_1:=\zeta_1(h_0)>0$. Let $\vp$ be the positive eigenfunction of (\ref{2.1}) corresponding to $\zeta_1$. Then $\vp'(h_0)<0$. There exists $0<x_0<h_0$ such that
 \bes
 \vp'(x)<0\ \ \ {\rm for}\ \ x\in[x_0,h_0].
  \lbl{2.8}\ees
We know that $\vp(x)>0$ on $[0,x_0]$ due to $b>0$. Hence, there exists a constant $m>0$ such that
 \bes
 \vp'(x)\le m\vp(x) \ \ \ {\rm for}\ \ 0\le x\le x_0.
 \lbl{2.9}\ees
 Let $0<\sigma\ll 1$ and $M>0$ be constants, which will be determined later. Set
 \bess
 \dd s(t)=1+2\sigma-\sigma {\rm e}^{-\sigma t}, \ \
 w(t,x)=M e^{-\sigma t}\vp\left(x/s(t)\right), \ \ \ t\geq 0,\ \ 0\leq x\leq h_0s(t).
 \eess
It is easy to see that $\vp$ also meet  (\ref{2.3}) and  $f(w)\le f'(0)w$. Denote $y=x/s(t)$, the direct calculations yield, for all $t>0$ and $0<x<h_0s(t)$,
\bes
&\ &w_t-w_{xx}+\beta w_x-f(w)\nonumber\\[.5mm]
&=&w\left(-\sigma- \frac{y\varphi'(y)s'(t)}{\varphi(y)s(t)}-\frac{\varphi''(y)}{\varphi(y)s^2(t)}+
\frac{\beta\varphi'(y)}{\varphi(y)s(t)}\right)-f(w) \nonumber\\[1mm]
 &\geq& w\left(-\sigma- \frac{y\varphi'(y)\sigma^2}{\varphi(y)s(t)}{\rm e}^{-\sigma t}+\frac{\zeta_1+f'(0)}{s^2(t)}
 +\frac{\beta\varphi'(y)\bbb(s(t)-1\bbb)}{\varphi(y)s^2(t)}-f'(0)\right).
 \label{2.10}\ees
We first estimate the term $\frac{\beta\varphi'(y)\bbb(s(t)-1\bbb)}{\varphi(y)s^2(t)}$. In view of (\ref{2.8}), $\beta\le 0$ and $s(t)>1$  we have
 \bes
 \frac{\beta\varphi'(y)\bbb(s(t)-1\bbb)}{\varphi(y)s^2(t)} \ge 0 \ \ \ {\rm for} \ x_0\le y<h_0.
  \lbl{2.11}\ees
By virtue of (\ref{2.9}), it follows that
  \bes
  \left|\frac{\beta\varphi'(y)\bbb(s(t)-1\bbb)}{\varphi(y)s^2(t)}\right|\le
   \left|\frac{\beta m\bbb(s(t)-1\bbb)}{s^2(t)}\right|\le 2\sigma|\beta|m \ \ \ {\rm for} \ 0\le y\le x_0.
  \lbl{2.12}\ees
In view of (\ref{2.3}), (\ref{2.11}),  (\ref{2.12}) and $\zeta_1>0$, it follows from (\ref{2.10}) that, for all $t>0$ and $0<x<h_0s(t)$,
 \bes
w_t-w_{xx}+\beta w_x-f(w)
\geq w\left(-\sigma-k\sigma^2+\frac{\zeta_1-4\sigma(1+\sigma)f'(0)}
{(1+2\sigma)^2}-2\sigma|\beta|m\right)>0
 \label{2.13}\ees
provided $0<\sigma\ll1$.

It is easy to see that $w(t,h_0g(t))=M{\rm e}^{-\sigma t}\varphi(h_0)=0$, and $B[w](t,0)=0$ if either $a=0$ or $b=0$. When $a>0$ and $b>0$, $B[\varphi](t,0)=0$ implies $a\varphi(0)=b\varphi'(0)$ and $\varphi'(0)>0$, so $B[w](t,0)=bM{\rm e}^{-\sigma  t}\varphi'(0)(1-1/s(t))>0$ on account of $s(t)>1$. In a ward,
\begin{equation}\label{2.14}
B[w](t,0)\geq0,\ \ \ w(t,h_0s(t))=0,\ \ \forall\  t\geq0.
\end{equation}

In addition, similar to Step 1 in the proof of Lemma \ref{lb.5}, we can select $M\gg 1$ and find $\mu_0=\mu_0(M)>0$ such that
\bes\label{2.15}
&u_0(x)\leq M\varphi\left( x/{(1+\sigma)}\right)=w(0,x), \ \ \forall\ 0\leq x\leq h_0,&\\[.1mm]
 &h_0s'(t)=h_0\sigma^2 {\rm e}^{-\sigma t}\geq -\mu w_x(t,h_0s(t)), \ \ \forall\ 0<\mu\leq\mu_0, \ t\geq0.&\label{2.16}
\ees

Combining (\ref{2.13})-(\ref{2.16}), we can apply Lemma \ref{l2.1} to achieve that
$h(t)\leq h_0s(t)$ for all $t\geq0$. Thus $h_\infty\leq h_0(1+2\sigma)$ for all $0<\mu\leq \mu_0$.
 \ \ \ \fbox{}

\begin{lemma}\lbl{lb.7}{\hspace{-1.8mm}\bf .} Let $|\beta|<c_0$ and $(u,h)$ be the solution of $(P)$. If $h_0<\pi/\sqrt{c^2_0-\beta^2}$, and either $\beta\leq0$, or $\beta>0$ and $a\geq b\beta/2$, then the conclusions of Theorem $\ref{th2.5}$ hold.
 \end{lemma}

\noindent{\bf Proof.} The proof is essentially the same to that of Lemma \ref{lb.5}, however the difference is the choice of auxiliary function. Thus we merely  show the following
  \vspace{-1mm}\begin{quote}{\bf Claim}: for any fixed $u_0\in \mathscr{X}(h_0)$, there exists $\mu_0>0$ depending on $\beta,h_0,f'(0)$ and $u_0(x)$ such that $h_\infty<\ty$ for $\mu\leq\mu_0$,\vspace{-1mm}\end{quote}
and the rest argument will not be duplicated here.

Let $0<\nu<1$ and $C>0$ be constants to be determined later. Set
 \bess
 p(t)&=&h_0(1+\nu-\frac{\nu}{2}e^{-\nu t}) \ \ \;{\rm for}\, \ t\geq0,\\
 z(t,x)&=&C e^{-\nu t}e^{-\frac\beta 2(p(t)-x)}\cos\frac{\pi x}{2p(t)} \ \ \ {\rm for}\, \ t\geq0, \ 0\leq x\leq p(t).\eess
Evidently, $z(t,p(t))=0$ and $p(0)=h_0(1+\frac{\nu}{2})>h_0$. Utilizing the assumptions on $\beta$ we can verify that $ B[z](t,0)\geq0$. Note that $h_0<\pi/\sqrt{c^2_0-\beta^2}$, i.e., $c_0^2-\beta^2<\pi^2/h_0^2$, and $f'(0)=4c_0^2$, straightforward computations generate that, for $t\geq0$ and $0\leq x\leq p(t)$,
\begin{eqnarray*}
\begin{array}{lll}
\displaystyle z_t- z_{xx}+\beta z_x-f( z)&= & \displaystyle z\left(-\nu-\frac\beta2p'(t)+\frac{\beta^2}4+\frac{\pi^2}{4p^2(t)}\right)
-f(z)\\[3mm]
&&+\displaystyle\frac{\pi x p'(t)}{2p^2(t)}C e^{-\nu t}e^{-\frac\beta2(p(t)-x)}\sin\frac{\pi x}{2p(t)}\\ [3mm]
&\geq&\displaystyle z\left(-\nu-\frac\beta4h_0\nu^2+\frac{\beta^2}4
+\frac{\pi^2}{4h_0^2(1+\nu)^2}-f'(0)\right)\geq0
\end{array}
\end{eqnarray*}
provided $0<\nu\ll 1$. Moreover, we can select $C\gg 1$ and find an $\mu_0=\mu_0(C)>0$ such that $u_0(x)\leq z(0,x)$ for all $0\leq x\leq h_0$, and
\bess
-\mu z_x(t,p(t))=\frac{\mu\pi C}{2p(t)}e^{-\nu t}\leq\frac{\mu\pi C}{2h_0}e^{-\nu t}\leq\frac{\nu^2h_0}{2}e^{-\nu t}=p'(t), \ \ \forall\ 0<\mu\leq\mu_0, \ t\geq0.
\eess

Similar to the above, we can employ Lemma \ref{l2.1} to derive
$h(t)\leq h_0s(t)$, and hence $h_\infty\leq h_0(1+2\sigma)$ for all $0<\mu\leq \mu_0$. \ \ \ \fbox{}

In what follows we present the sharp criteria on initial value, which separates vanishing and spreading.

\begin{theorem}\label{th2.6}\hspace{-1mm}{\rm(}Threshold on initial value{\rm)}. Assume that $|\beta|<c_0$ and $(u,h)$ is a solution of $(P)$ with $u_0=\lambda \psi$, where $\lambda>0$ and $\psi\in\mathscr{X}(h_0)$. Then there exists $\lambda^*\in[0,\infty]$ dependent on $h_0, f,\psi$ so that spreading happens when $\lambda>\lambda^*$, and vanishing happens when $0<\lambda\leq\lambda^*$ provided $\lambda^*>0$.

Furthermore,  $\lambda^*=0$ if $h_0\ge\ell^*$, and $\lm^*>0$ if the conditions of Theorem $\ref{th2.5}$ hold.
\end{theorem}

\noindent{\bf Proof.}
Define
$$\sum=\big\{\lambda_0>0: {\rm vanishing\ happens\ for\ all}\ \lambda\in(0,\lambda_0]\big\}$$
and $\lambda^*=\sup\sum$. If $h_0\ge\ell^*$, we have $\sum=\emptyset$ by Theorem \ref{th2.3}, and set $\lambda^*=0$. When $\sum=(0,\infty)$, then $\lambda^*=\infty$, which implies that vanishing happens no matter how large $\lambda$ is. In case $\lambda^*=\infty$ (let us point out that this happens in particular when $\beta=0$, $\liminf_{s\to\infty}\frac{-f(s)}{s^{8}}\gg1$ and $a=0$, refer to \cite[Proposition 5.4]{DLou} for details), there is nothing left to show.

In the following we suppose $0<\lambda^*<\infty$. According to the definition of $\lambda^*$ and spreading-vanishing dichotomy, we can find a sequence $\lambda_i$ decreasing to $\lambda^*$ so that spreading happens when $\lambda=\lambda_i,\ i=1,2,\cdots$. For any given $\lambda>\lambda^*$, we can select some $i\geq1$ so that $\lambda>\lambda_i$. Denote by $(u_i,h_i)$ the solution of problem $(P)$ with $u_0=\lambda_i\psi$. Then, in terms of the comparison principle, we know that $[0,h_i(t)]\subset[0,h(t)]$ and $u_i(t,x)\leq u(t,x)$ for all $t>0$ and $0\leq x\leq h_i(t)$. Therefore, spreading happens for such $\lambda$. We shall demonstrate that vanishing happens when $\lambda=\lambda^*$. Suppose on the contrary that spreading happen for $\lambda=\lambda^*$. Then there exists $t_0>0$ such that $h(t_0)>\ell^*+1$. Utilizing the continuous dependence of the solution for problem $(P)$ on its initial values, we can choose $\varepsilon>0$ sufficiently small such that the solution of $(P)$ with $u_0=(\lambda-\varepsilon)\psi$, denoted by $(u_\varepsilon,h_\varepsilon)$, satisfies
   $$h_\varepsilon(t_0)>\ell^*.$$
In view of Theorem \ref{th2.3} we see that spreading happens for $(u_\varepsilon,h_\varepsilon)$. This leads to a contradiction with the definition of $\lambda^*$.

At last, we demonstrate that $\lm^*>0$ if the conditions of Theorem $\ref{th2.5}$ hold. In fact, it suffices to show
 \vspace{-1mm}\begin{quote}{\bf Claim}:  If the conditions of Theorem \ref{th2.5} hold, then there exists $\lambda_0>0$ depending on $\beta,h_0,\gamma_1(h_0)$ and $\mu$, such that $h_\infty<\infty$ for $0<\lambda\leq \lambda_0$.\vspace{-1mm}\end{quote}
In the same way to Theorem \ref{th2.5}, the argument for this claim will divided into three cases:

(i) $h_0<\ell^*$ and $b=0$;

(ii) $h_0<\ell^*$, and $b>0$, $\beta\le 0$;

  (iii) $h_0<\pi/\sqrt{c^2_0-\beta^2}$, and either $\beta\leq0$, or $\beta>0$ and $a\geq b\beta/2$.\\
Because the proofs for case (i), case (ii) and case (iii) are similar to Step 1 of Lemma \ref{lb.5}, Lemma \ref{lb.6} and Lemma \ref{lb.7}, respectively, we only provide a sketch for argument of case (i).

 Let $\gamma_1$ and $\phi$ be as above, and $0<\delta,\varepsilon<1$ be constants to be determined later. Set
 \bess
 \dd g(t)=1+\delta-\frac\delta 2{\rm e}^{-\delta t}, \ \
 v(t,x)=\varepsilon e^{-\delta t}e^{-\frac\beta2(h_0g(t)-x)}\phi\left(x/g(t)\right), \ \ \ t\geq 0,\ \ 0\leq x\leq h_0g(t).
 \eess
Similar to Step 1 of Lemma \ref{lb.5}, we can see that (\ref{2.4}) and (\ref{2.5}) still hold provided $0<\delta\ll1$. Fix the $\delta$, we can choose $0<\varepsilon\ll 1$ such that, for all $t\geq0$,
  $$-\mu v_x(t,h_0g(t))=-\frac{\mu\varepsilon\phi'(h_0)}{g(t)}{\rm e}^{-\delta t}\leq -\mu\varepsilon\phi'(h_0){\rm e}^{-\delta t}\leq h_0\delta^2 {\rm e}^{-\delta t}= h_0g'(t).$$
Further fix the $\varepsilon$. By the regularities of $\psi(x)$ and $\phi(x)$, we can take $0<\lambda_0\ll 1$ such that, for all $0<\lambda\leq\lambda_0$,
  \[u_0(x)=\lambda \psi(x)\leq \varepsilon e^{-\frac\beta2(h_0(1+\frac\delta2)-x)}\phi\left( 2x/{(2+\delta)}\right)=v(0,x), \ \ \forall\ 0\leq x\leq h_0.\]

In summary,
\begin{eqnarray*}
\left\{\begin{array}{ll}
v_t-v_{xx}+\beta v_x-f(v)\geq0, & t\geq 0,\ 0\leq x\leq h_0g(t),\\[1mm]
B[v](t,0)\geq 0, \ \ v(t,h_0g(t))=0, \ \ &t\geq 0,\\[1mm]
h_0g'(t)\geq -\mu v_x(t,h_0g(t)),&t\geq0,\\[1mm]
v(0,x)\geq  u_0(x), & 0\leq x\leq h_0.
\end{array}\right.
\end{eqnarray*}
Now we can easily obtain the required result in the same manner to Step 1 of Lemma \ref{lb.5}. This completes the proof of this theorem.
\ \ \ \fbox{}

\vspace{.4cm}
Based upon Theorems \ref{th2.3} and \ref{th2.5}, we derive the following sharp criteria for spreading and vanishing.

\begin{corollary}\hspace{-2mm}{\bf .}\lbl{c2.1} Assume that $(u,h)$ be the solution of $(P)$, and that $b=0$ and $|\beta|<c_0$, or $b>0$ and $-c_0<\beta\le 0$. We have

{\rm(i)}\, if $h_0\ge\ell^*$, then $h_\ty=\ty$ for any $u_0\in \mathscr{X}(h_0)$ and $\mu>0$;

{\rm(ii)} if $h_0<\ell^*$,
then the following assertions hold:

\vskip 4pt {\rm(iia)} For any given $u_0\in \mathscr{X}(h_0)$, there exists $\mu^*>0$
such that spreading happens for $\mu>\mu^*$, while vanishing happens for $\mu\leq\mu^*$.

\vskip 4pt
{\rm(iib)} Fix $\mu>0$ and take $u_0=\lambda \psi$ with $\lambda>0$ and $\psi\in\mathscr{X}(h_0)$. Then there exists $\lambda^*\in(0,\infty]$
so that vanishing happens when $0<\lambda\leq\lambda^*$, and spreading happens when $\lambda>\lambda^*$.
\end{corollary}

\section{Uniform convergence for $u$ and sharp estimates of $h(t)$ and $h'(t)$} \label{se3}
\setcounter{equation}{0}

Throughout this section we assume that $|\beta|<c_0$ and $(u,h)$ is a solution of $(P)$ for which spreading happens. Consider the following elliptic problem:
\begin{equation}\label{a3.1}
\left\{\begin{array}{l}
q''-(c-\beta)q'+f(q)=0,\ \ 0<z<\infty,\\[1mm]
q(0)=0,\ \ q(\infty)=1,\\[1mm]
 q(z)>0,\ \ 0<z<\infty.
\end{array}
\right.
\end{equation}

\begin{proposition}\label{p3.1}{\hspace{-1.8mm}\bf .} For any given $\mu>0$, there exist a unique $\tilde c_\beta=\tilde c_\beta(\mu)\in(0,c_0+\beta)$ and a unique solution $\tilde q_\beta(z)$ to $(\ref{a3.1})$ with $c=\tilde c_\beta$ such that $\tilde q_\beta'(0)=\tilde c_\beta/\mu$. Moreover, $\tilde c_\beta$ is increasing in $\beta$ and $\ds\lim_{\beta\to-c_0}\tilde c_\beta=0$.
\end{proposition}

\noindent{\bf Proof.}
Since the proof is essentially identical to that of \cite[Lemma 3.4]{GLZ}, we only present a sketch here for the reader's convenience.

For any given $c<c_0+\beta$, problem (\ref{a3.1}) admits a unique solution, denoted by $q_\beta(z;c-\beta)$, which satisfies $q_\beta'(z;c-\beta)>0$ for $0\le z<\infty$. Denote $P(c-\beta)=\mu q_\beta'(0;c-\beta)$, then $P(c-\beta)$ is strictly decreasing in $c\in(-\infty,c_0+\beta)$, and
  $$(P(c-\beta)-c)|_{c=0}=P(-\beta)>0, \ \ \ \lim_{c\to c_0+\beta-0}(P(c-\beta)-c)=-c_0-\beta<0.$$
Thus the equation $P(c-\beta)=c$ admits a unique solution $c=\tilde c_\beta\in(0,c_0+\beta)$. Let $\tilde q_\beta(z)=q_\beta(z;\tilde c_\beta-\beta)$. Then we have
$$\tilde c_\beta=P(\tilde c_\beta-\beta)=\mu\tilde q_\beta'(0;\tilde c_\beta-\beta).$$

Differentiating $\tilde c_\beta=P(\tilde c_\beta-\beta)$ and utilizing the fact $P'(c)<0$ for $c<c_0$ we obtain
$$\frac{d\tilde c_\beta}{d\beta}=\frac{-P'(\tilde c_\beta-\beta)}{1-P'(\tilde c_\beta-\beta)}>0.$$
Hence $\tilde c_\beta$ is increasing in $\beta$. It is easy to find $\ds\lim_{\beta\to-c_0}\tilde c_\beta=0$ from the above range of $\tilde c_\beta$.
The proof is finished.
\ \ \ \fbox{}

\begin{remark}\label{k3.1}{\hspace{-1.8mm}\bf .} The transformation of independent variables  enables problem  $(\ref{a3.1})$ to become
\begin{eqnarray*}
\left\{\begin{array}{l}
q''+(c-\beta)q'+f(q)=0,\ \ -\infty<z<0,\\[1mm]
q(0)=0,\ \ q(-\infty)=1,\\[1mm]
 q(z)>0,\ \ -\infty<z<0.
\end{array}
\right.
\end{eqnarray*}
Therefore, combining Lemma $3.4$ of {\rm \cite{GLZ}} and Proposition $\ref{p3.1}$, one can know that, for any $\mu>0$ and $-c_0<\beta<\infty$, there exist a unique $\tilde c_\beta\in(0,c_0+\beta)$ and a unique solution $\tilde q_\beta(z)$ to $(\ref{a3.1})$ with $c=\tilde c_\beta$ such that $\tilde q_\beta'(0)=\tilde c_\beta/\mu$. Besides, Lemma $3.4$ of {\rm \cite{GLZ}}  illustrates that there exists a unique $\beta^*>c_0$ such that $\tilde c_\beta-\beta+c_0>0$ $(resp.\ =0,\ <0)$ when $\beta<\beta^*$ $(resp.\ \beta=\beta^*,\ \beta>\beta^*)$.
\end{remark}
\begin{remark}\label{k3.2}{\hspace{-1.8mm}\bf .}\ From the argument of Proposition \ref{p3.1} we see that $\tilde c_\beta=\tilde c_\beta(\mu)$ is the unique solution of
 $$c=\mu q_\beta'(0;c-\beta). $$
Fix $\beta$ with $|\beta|<c_0$. Then $\tilde c_\beta(\mu)\to 0$ as $\mu\to0$ since $q_\beta'(0;c-\beta)$ is bounded in $c\in[0, c_0+\beta]$.
\end{remark}

The function $\tilde q_\beta(z)$ obtained in Remark \ref{k3.1} is called a semi-wave with speed $\tilde c_\beta$ since $u(t,x)=\tilde q_\beta(\tilde c_\beta t-x)$ satisfies
\begin{eqnarray*}
\left\{\begin{array}{l}
u_t-u_{xx}+\beta u_x=f(u), \ \ \ t\in \mathbb{R},\ \ x<\tilde c_\beta t,\\[1mm]
u(t,-\infty)=1,\ u(t,\tilde c_\beta t)=0,\ -\mu u_x(t,\tilde c_\beta t)=\tilde c_\beta.
\end{array}\right.
\end{eqnarray*}
In what follows we introduce a wave of finite length which will be used to construct lower solutions of problem $(P)$. Employing a similar phase plane analysis as in \cite{DLou} we can achieve the following proposition.

\begin{proposition}{\hspace{-1.8mm}\bf .} For every $c\in(0,\tilde c_\beta)$, there exists a unique pair $(q_c(z),z_c)$ satisfying
 \begin{equation}\label{a3.2}
\left\{\begin{array}{ll}
{q_c} ''-(c-\beta) {q_c} '+f(q_c)=0,\ \ \ z\in[0,z_c],\\[1mm]
q_c(0)=0,\ \ \mu {q_c} '(0)=\tilde c_\beta,\ \ q_c'(z_c)=0,\\[1mm]
q_c(z)>0\ \ \ {\rm in}\, \ (0,z_c].
\end{array}
\right.
\end{equation}
Moreover, $\ds\lim_{c\nearrow \tilde c_\beta}z_c=\infty$ and $\ds\lim_{c\nearrow \tilde c_\beta}\|q_c-\tilde q_\beta\|_{L^\infty([0,z_c])}=0$.
\end{proposition}

The next conclusion plays an important role in arriving at the bound of $h(t)-\tilde c_\beta t$ and the uniform convergence of $u(t,x)$.

\begin{lemma}\label{le4}{\hspace{-1.8mm}\bf .}
Assume that $(u,h)$ is a solution of $(P)$ for which spreading happens. Then for any $c\in(0,\tilde c_\beta)$, $\sigma\in(0,-f'(1))$, there exist positive constants $T$, $K$, $C$, $\bar c\in(0,c/2)$ and $\bar\sigma=\bar\sigma(c)\in(0,\sigma)$  such that, for $t\geq T$,

{\rm(i)}
 $h(t)\geq ct$;

{\rm(ii)}
$u(t,x)\leq 1+Ke^{-\sigma t}$  for all  $0\leq x\leq h(t)$;

{\rm(iii)} $u(t,x)\geq 1-Ce^{-\bar\sigma t}$ for all $\frac{(c-\bar c)t}2\leq x\leq \frac{(c+\bar c)t}2$.
\end{lemma}

\noindent{\bf Proof.} (i) Fix $\hat c\in(c,\tilde c_\beta)$ and define
 \bess
 &g(t)=z_{\hat c}+\hat ct+x_0\ \ {\rm for }\ t>0,\\
 &w(t,x)=q_{\hat c}(g(t)-x)\ \ {\rm for }\ t>0,\ g(t)-z_{\hat c}\leq x\leq g(t),&\eess
where $(q_{\hat c}(z),z_{\hat c})$ satisfies (\ref{a3.2}) with $c=\hat c$, and $x_0$ is to be determined later. Since spreading happens, we can select $T_1>0$ and $x_0>0$ such that $h(T_1)\geq g(0),$ $$u(t,g(t)-z_{\hat c})>w(t,g(t)-z_{\hat c})\ \ {\rm for }\, \ t\geq T_1,\ \ \ u(T_1,x)>w(0,x)\ \ {\rm in}\, \ [g(0)-z_{\hat c},g(0)].$$
It is easy to verify that $(w(t-T_1,x),g(t-T_1))$ is a lower solution of $(P)$ for $t\geq T_1$ and $g(t)-z_{\hat c}\leq x\leq g(t)$. Therefore, by Lemma \ref{l2.4} and the remark behind it we can deduce that, for some $T_2>T_1$,
$$h(t)\geq g(t-T_1)>\hat c(t-T_1)>ct\ \ {\rm for} \ t\geq T_2.$$

(ii) Consider the equation $\eta'(t)=f(\eta)$ with $\eta(0)=\|u_0\|_{L^\infty}+1$. It is readily seen that $\eta$ is an upper solution of $(P)$. Therefore, $u(t,x)\leq \eta(t)$ for all $t>0$. Since $f(u)<0$ for $u>1$, $\eta(t)$ is a decreasing function satisfying
\begin{equation}\label{a3.3}
\lim_{t\to\infty}\eta(t)=1.
\end{equation}
Noting that $0<\sigma<-f'(1)$,  we can choose some $0<\varsigma<1$ so that
\begin{eqnarray}\label{a3.4}
\left\{\begin{array}{ll}
\sigma\leq-f'(u)\ & {\rm for}\ \;1-\varsigma\leq u\leq 1+\varsigma,\\[1mm]
f(u)\geq0\ & {\rm for}\ \; 1-\varsigma\leq u\leq1.
\end{array}\right.
\end{eqnarray}
It follows from (\ref{a3.3}) that there is a positive number $T_3$ such that $\eta(t)<1+\varsigma$ for all $t\geq T_3$. And then $\eta'(t)=f(\eta)\leq \sigma(1-\eta)$ for all $t\geq T_3$. Now it is not difficult to obtain that
  $$u(t,x)\leq \eta(t)\leq1+\varsigma e^{\sigma T_3} e^{-\sigma t}\  \ {\rm for \ } t\geq T_3,\ 0\leq x\leq h(t).$$

{(iii)} Define
\bess
 &g(t)\ds=-\frac c2 t,\ \ s(t)=h(t)-\frac c2 t\ \ \ {\rm for}\ \ t\geq0,&   \\[1mm]
 &w(t,y)=u(t,y+\frac c2 t) \ \ \ {\rm for}\ \ t\geq0\ \ {\rm and}\ \ y\in[g(t),s(t)].&
\eess
Then $ (w,g,s)$ solves the problem
\begin{equation}\label{a3.5}
\left\{
\begin{array}{ll}
 w_t- w_{yy}+(\beta-c/2) w_y=f( w),&t>0,\ g(t)<y<s(t),\\[1mm]
 B[w](t,g(t))=0,\ g'(t)=- c/2,&t>0,\\[1mm]
 w(t,s(t))=0,\ s'(t)=-\mu w_y(t,s(t))-c/2,\ \ &t>0,\\[1mm]
g(0)=0,\ \ s(0)=h_0,\ \  w(0,y)=u_0(y),&0\leq y\leq h_0.
\end{array}\right.
\end{equation}
Since $g(t)=-\frac c2 t\to-\infty$ and $s(t)=h(t)-\frac c2 t\geq \frac c2t\to\infty$ as $t\to\infty$, spreading happens for (\ref{a3.5}).  Let  $\underline s_0>0$ to be chosen later,  and $\underline w_0(x)\in C^2([-\underline s_0,\underline s_0])$ satisfying $$\underline w_0(-\underline s_0)=\underline w_0(-\underline s_0)=0,\ \underline w_0'(-\underline s_0)>0,\ \underline w_0'(\underline s_0)<0,\  \underline w_0(x)>0\ {\rm in}\ (-\underline s_0,\underline s_0).$$ Consider the following free boundary problem
\begin{equation}\label{a3.55}
\left\{
\begin{array}{ll}
 \underline w_t- \underline w_{yy}+\hat\beta \underline w_y=f(\underline  w),&t>0,\ \underline g(t)<y<\underline s(t),\\[1mm]
 \underline w(t,\underline g(t))=0,\ \underline g'(t)=- \underline\mu\underline w_y(t,\underline g(t)) ,&t>0,\\[1mm]
 \underline w(t,\underline s(t))=0,\ \underline s'(t)=-\underline \mu \underline w_y(t,\underline s(t)),\ \ &t>0,\\[1mm]
\underline g(0)=-\underline s_0,\ \ \underline s(0)=\underline s_0,\ \  \underline w(0,y)=\underline w_0(y),&-\underline s_0\leq y\leq \underline s_0,
\end{array}\right.
\end{equation}
where $\underline \mu>0$ and $\hat\beta=\beta-c/2\in(-c_0,c_0)$.  Taking advantage of a similar argument to Theorem \ref{th2.3}, one can find a large $\underline s_0$ independent of $\underline\mu$ and $\underline w_0$ such that spreading happens for the problem (\ref{a3.55}). For such fixed $\underline s_0$, combining Theorem B of \cite{KM} with Proposition \ref{p3.1}, it follows that there exist $\tilde c_{-\hat\beta}(\underline\mu)\in(0,c_0-\hat\beta)$, $\tilde c_{\hat\beta}(\underline\mu)\in(0,c_0+\hat\beta)$ such that $\dd\lim_{t\to\infty}\underline g'(t)=-\tilde c_{-\hat\beta}(\underline\mu),\ \,\lim_{t\to\infty}\underline s'(t)=\tilde c_{\hat\beta}(\underline\mu)$. According to Remark \ref{k3.2}, $\ds\lim_{\underline\mu\to 0}\tilde c_{-\hat\beta}(\underline\mu)=\lim_{\underline\mu\to 0}\tilde c_{\hat\beta}(\underline\mu)=0$. There exists a small $\underline\mu$ ($0<\underline\mu\le\mu$) independent of $\underline w_0$ such that $\tilde c_{-\hat\beta}(\underline\mu)$, $\tilde c_{\hat\beta}(\underline\mu)<c/2$. Hence, there exists $T_4>0$ large enough such that
  $$g(t+T_4)\leq \underline g(t), \ \, \underline s(t)\leq s(t+T_4)\ \,{\rm for\ all }\ t\geq 0.$$
Now we can select $\underline w_0(x)\leq w(T_4,x)$, then the solution $(\underline w,\underline g,\underline s)$ of problem (\ref{a3.55}) for which spreading happens is a lower solution of (\ref{a3.5}) for $t\geq T_4$.
By the same argument as \cite[Proposition 3.2]{KM} with some obvious modifications, we can still demonstrate that for any $0<c'<\min\{\tilde c_{-\hat\beta}(\underline\mu), \tilde c_{\hat\beta}(\underline\mu)\}$, there exist $\bar\sigma\in(0,\sigma)$ and $C'$, $T_5>0$ such that $\underline w(t,y)\geq 1-C'e^{-\bar\sigma t}$ for all $t>T_5$ and $-c't/2\leq y\leq c't/2$. It is not too difficult to show that there exists ${\bar c}\in(0,c')$, $C>0$ such that for some $T_6>T_4+T_5$,
  $$u(t,x)\geq 1-Ce^{-\bar\sigma t}\ \ \ {\rm for}\ t>T_6,\ \ \frac{(c-\bar c)t}2\leq x\leq \frac{(c+\bar c)t}2.$$

Denote $T=\max\{T_2,T_3,T_6\}$, and  $T$  can be used in place of $T_2$, $T_3$ and $T_6$. This completes the proof.
\ \ \ \fbox{}

Next, based on Lemma \ref{le4} we are going to construct two groups of upper and lower solutions for problem $(P)$, which play a crucial role in achieving the uniform convergence of $u(t,x)$.

Let $\tilde c_\beta$ and $\tilde q_\beta(z)$ be determined by Proposition \ref{p3.1}.  Let $c\in(0,\tilde c_\beta)$, $\sigma\in(0,-f'(1))$, the positive constants $T$, $\bar c$ and $\bar\sigma$ be obtained by Lemma \ref{le4}. For convenience of discussion, denote
  \begin{equation}\label{alr}c_l=\frac{c-\bar c}2,\ \ c_r=\frac{c+\bar c}2.\end{equation}
Taking advantage of the estimates (ii) and (iii) in Lemma \ref{le4}, one can directly calculate as in \cite{KM} (also see \cite{DMZ}) to demonstrate the following conclusion.

\begin{lemma}\label{le5}{\hspace{-1.8mm}\bf .} There exist positive constants $K_i, M_i$, $i=1,2$, $\kappa\gg 1$ and $T^*>T$ such that
\begin{eqnarray*}h(t)\leq \bar h(t)\ \;{\rm for}\ t\geq T^*,\ \;u(t,x)\leq\bar u(t,x)\ \; {\rm for}\ t\geq T^*,\ \bar g(t)\leq x\leq h(t);\\
h(t)\geq \underline h(t)\ \;{\rm for}\ t\geq T^*,\ \;u(t,x)\geq\underline u(t,x)\ \; {\rm for}\ t\geq T^*,\ \underline g(t)\leq x\leq \underline h(t),\end{eqnarray*}
where\bes
  &\bar g(t)= 0,\ \ \ \bar h(t)=\tilde c_\beta t-\kappa K_1 e^{-\sigma t}+M_1,\ \ \ \bar u(t,x)=(1+K_1e^{-\sigma  t})\tilde q_\beta(\bar h(t)-x),&\label{3.6}\\[1mm]
  &\underline g(t)=c_lt,\ \ \ \underline h(t)=\tilde c_\beta t+\kappa  K_2 e^{-\bar\sigma t}+M_2,\ \ \ \underline u(t,x)=(1- K_2e^{-\bar\sigma  t})\tilde q_\beta(\underline h(t)-x). &\label{3.7}\ees
  \end{lemma}

For the case $a>0$, other group of upper and lower solutions are constructed below. It follows from Lemma \ref{le4} that, for any $\sigma\in(0,-f'(1))$, there exist $T$, $K>0$ such that
\begin{equation}\label{a3.8}u(t,x)\leq 1+Ke^{-\sigma t}\ \ \ {\rm for}\ \, t\geq T,\ 0\leq x\leq h(t).\end{equation}
Let $\varsigma$ be given in (\ref{a3.4}). Then we can enlarge $T$ such that
$Ke^{-\sigma T}\leq\varsigma/2.$
We further choose $M>K$ satisfying
$Me^{-\sigma T}\leq\varsigma.$
For such a $\varsigma$, we define $x_\varsigma\in(0,\infty)$ and $V_\varsigma$ as follows:
  \begin{equation}\label{a3.9}\tilde v(x_\varsigma)=1-\varsigma,\ \;V_\varsigma=\min_{0\leq x\leq x_\varsigma}\tilde v'(x),\end{equation}
where $\tilde v(x)$ is the unique strictly increasing solution of (\ref{1.4}) on the half-line with $\tilde v(\infty)=1$. There exists $x^*>0$ so that
\begin{equation}\label{a3.10}(1+Me^{-\sigma T})\tilde v(x^*)\geq 1+Ke^{-\sigma T}.\end{equation}
Define
\bes
 \bar l(t)=\kappa M(e^{-\sigma t}- e^{-\sigma T})-x^*,\ \
 \bar u(t,x)=(1+Me^{-\sigma t})\tilde v(x-\bar l(t)),
 \label{3.11}\ees
where $\kappa$ is a positive constant to be determined.

Notice that $\bar l(t)<0$ for $t\geq T$ and $B[\tilde v](0)=0$. By the property of $\tilde v(x)$, it is easy to see  that $\bar u(t,h(t))>0= u(t,h(t))$ for $t\geq T$ and $B[\bar u](t,0)=B[\tilde v](-\bar l(t))\geq0$ for $t\geq T$. According to (\ref{a3.8}) and (\ref{a3.10}), we can get, for $0\leq x\leq h(T)$,
  \[\bar u(T,x)=(1+Me^{-\sigma T})\tilde v(x+x^*)
   \geq(1+Me^{-\sigma T})\tilde v(x^*)\geq 1+Ke^{-\sigma T}\geq  u(T,x).\]

Set $z=x-\bar l(t)$. A routine calculation shows that, for $t\geq T$ and $\bar l(t)\leq x\leq h(t)$,
\begin{eqnarray}
&&\bar u_t-\bar u_{xx}+\beta \bar u_x-f(\bar u)\nonumber  \\[1mm]
\nonumber&=&-\sigma M e^{-\sigma t}\tilde v(z)-(1+Me^{-\sigma t})\bar l'(t)\tilde v'(z)-(1+Me^{-\sigma t})\tilde v''(z)\\[1mm]
\nonumber&&+\beta (1+Me^{-\sigma t})\tilde v'(z)-f((1+Me^{-\sigma t})\tilde v(z))\\[1mm]
\nonumber           &=&-\sigma M e^{-\sigma t}\tilde v(z)+\kappa M\sigma e^{-\sigma t}(1+Me^{-\sigma t})\tilde v'(z)\\[1mm]
 \nonumber          &&+(1+Me^{-\sigma t})f(\tilde v(z))-f((1+Me^{-\sigma t})\tilde v(z))\\[1mm]
\nonumber           &=&-\sigma M e^{-\sigma t}\tilde v(z)+\kappa M\sigma e^{-\sigma t}(1+Me^{-\sigma t})\tilde v'(z)+Me^{-\sigma t}f(\tilde v(z))\\[1mm]
         &&-f'(\tilde v(z)+\theta Me^{-\sigma t}\tilde v(z))Me^{-\sigma t}\tilde v(z),\label{a3.11}
\end{eqnarray}
where we have applied the  mean value theorem with $0<\theta<1$. By virtue of (\ref{a3.9}) and $Me^{-\sigma t}\leq \varsigma$ for all $t\geq T$, we can deduce that, for $z=x-\bar l(t)$,
$$1-\varsigma\leq \tilde v(z)+\theta Me^{-\sigma t}\tilde v(z)\leq \tilde v(z)+\theta\varsigma \tilde v(z)<1+\varsigma.$$
Thus it follows from (\ref{a3.4}) and (\ref{a3.11}) that for $z\geq x_\varsigma$ and $t\geq T$,
$$\bar u_t-\bar u_{xx}+\beta \bar u_x-f(\bar u)\geq Me^{-\sigma t}\tilde v(z)\big[-\sigma-f'(\tilde v(z)+\theta Me^{-\sigma t}\tilde v(z))\big]\geq0.$$
On the other hand, from (\ref{a3.11}) we have, for $0\leq z\leq x_\varsigma$,
\begin{eqnarray*}
\bar u_t-\bar u_{xx}+\beta \bar u_x-f(\bar u)
           &\geq&\kappa M\sigma e^{-\sigma t}(1+Me^{-\sigma t})V_\varsigma -Me^{-\sigma t}\left(\sigma+\max_{0\leq s\leq 1+\varsigma}f'(s)\right)\\[1mm]
           &\geq& Me^{-\sigma t}\left(\kappa \sigma V_\varsigma-\sigma-\max_{0\leq s\leq 1+\varsigma}f'(s)\right).
\end{eqnarray*}
Hence we can choose  sufficiently large $\kappa>0$ such that $\bar u_t-\bar u_{xx}+\beta \bar u_x-f(\bar u)\geq0$ for $t\geq T$ and for $0\leq z\leq x_\varsigma$.

Note that $\bar l(t)<0$ for $t\ge T$. Summarizing the above discussion we obtain
\begin{eqnarray*}
\left\{\begin{array}{ll}
\bar u_t-\bar u_{xx}+\beta \bar u_x-f(\bar u)\geq0, & t\geq T,\ 0\leq x\leq h(t),\\[1mm]
B[\bar u](t,0)\geq 0,\ \bar u(t,h(t))\geq u(t,h(t)),\ \ &t\geq T,\\[1mm]
\bar u(T,x)\geq  u(T,x), & 0\leq x\leq h(T).
\end{array}\right.
\end{eqnarray*}
Thus the comparison principle enables us to conclude that

\begin{lemma}\label{le6}{\hspace{-1.8mm}\bf .} If $a>0$, then $u(t,x)$ satisfies
$u(t,x)\leq \bar u(t,x)\ \;{\rm for}\ t\geq T,\ 0\leq x\leq h(t)$, where $\bar  u(t,x)$ is given by $(\ref{3.11})$ and $\kappa>0$ is suitable large.
\end{lemma}

Fix $c\in(0,\tilde c_\beta)$ and $\sigma\in(0,-f'(1))$. According to Lemma \ref{le4}, there exist $\bar\sigma\in(0,\sigma)$, $\bar c\in(0,c/2)$ and $T$, $C>0$ such that
\begin{equation}\label{a3.12}u(t,x)\geq 1-Ce^{-\bar\sigma t}\ \ \ {\rm for}\ t\geq T,\ c_lt\leq x\leq c_rt,\end{equation}
where $c_l$ and $c_r$ are defined by (\ref{alr}).
 Let $\bar\varsigma$ satisfy (\ref{a3.4}) with $\sigma$ replaced by $\bar\sigma$.
By enlarging $T$ we can suppose that
 \begin{equation}\label{a3.13}Ce^{-\bar\sigma T}<\bar\varsigma/2.\end{equation}
Define
 \bes
\underline r(t)=c_rt, \ \ \underline l(t)=c_lT-\kappa C(e^{-\bar\sigma t}- e^{-\bar\sigma T}),\ \
\underline u(t,x)=(1-Ce^{-\bar\sigma t})\tilde v_0(x-\underline l(t)),\label{3.15}
 \ees
where $\tilde v_0(x)$ is the strictly increasing solution of problem (\ref{1.4}) with $b=0$ and $\tilde v_0(\infty)=1$.
Now we define $x_{\bar\varsigma/2}\in(0,\infty)$ and $V_{\bar\varsigma/2}$ as follows:
 \begin{equation}\label{a3.14}\tilde v_0(x_{\bar\varsigma/2})=1-\bar\varsigma/2,\ \;V_{\bar\varsigma/2}=\min_{0\leq x\leq x_{\bar\varsigma/2}}\tilde v_0'(x)>0.\end{equation}

Note that $\underline l(t)>0$ for all $t\geq T$ and  $B[\tilde v_0](0)=\tilde v_0(0)=0$ when $b=0$. Then it is easy to see that $\underline u(t,\underline l(t))=0\leq  u(t,\underline l(t)$ for $t\geq T$. By (\ref{a3.12}), we have $\underline u(t,\underline r(t))\leq (1-Ce^{-\bar\sigma t})\leq  u(t,\underline r(t))$ for $t\geq T$ and $\underline u(T,x)\leq  u(T,x)$ for $\underline l(T)\leq x\leq \underline r(T)$.

Set $y=x-\underline l(t)$. Using a similar calculation to (\ref{a3.11}) generates that
\begin{eqnarray}
&&\underline u_t-\underline u_{xx}+\beta \underline u_x-f(\underline u)\nonumber  \\[1mm]
\nonumber&=&\bar\sigma C e^{-\bar\sigma t}\tilde v_0(y)-(1-Ce^{-\bar\sigma t})\underline l'(t)\tilde v_0'(y)-(1-Ce^{-\bar\sigma t})\tilde v_0''(y)\\[1mm]
\nonumber&&+\beta (1-Ce^{-\bar\sigma t})\tilde v_0'(y)-f((1-Ce^{-\bar\sigma t})\tilde v_0(y))\\[1mm]
\nonumber           &=&\bar\sigma C e^{-\bar\sigma t}\tilde v_0(y)-\kappa C\bar\sigma e^{-\bar\sigma t}(1-Ce^{-\bar\sigma t})\tilde v_0'(y)\\[1mm]
 \nonumber          &&+(1-Ce^{-\bar\sigma t})f(\tilde v_0(y))-f((1-Ce^{-\bar\sigma t})\tilde v_0(y))\\[1mm]
\nonumber           &=&\bar\sigma C e^{-\bar\sigma t}\tilde v_0(y)-\kappa C\bar\sigma e^{-\bar\sigma t}(1-Ce^{-\bar\sigma t})\tilde v_0'(y)-Ce^{-\bar\sigma t}f(\tilde v_0(y))\\[1mm]
         &&+f'(\tilde v_0(y)-\varrho Ce^{-\bar\sigma t}\tilde v_0(y))Ce^{-\bar\sigma t}\tilde v_0(y)\label{a3.15}
\end{eqnarray}
for $t\geq T$ and $\underline l(t)\leq x\leq \underline r(t)$, where $0<\varrho<1$. We first discuss the case $y\geq x_{\bar\varsigma/2}$. According to (\ref{a3.13}) and (\ref{a3.14}), it holds that
$$1\geq \tilde v_0(y)-\varrho Ce^{-\bar\sigma t}\tilde v_0(y)\geq (1- Ce^{-\bar\sigma t})\tilde v_0(y)\geq(1-\bar\varsigma/2)^2\geq 1-\bar\varsigma.$$
Therefore, it follows from (\ref{a3.4}) with $\sigma$ replaced by $\bar\sigma$ that $f'(\tilde v_0(y)-\varrho Ce^{-\bar\sigma t}\tilde v_0(y))+\bar\sigma\leq0$. Consequently, we can know from (\ref{a3.15}) that, for $t\geq T$  and $y\geq x_{\bar\varsigma/2}$,
$$\underline u_t-\underline u_{xx}+\beta \underline u_x-f(\underline u)\leq Ce^{-\bar\sigma t}\tilde v_0(y)[f'(\tilde v_0(y)-\varrho Ce^{-\bar\sigma t}\tilde v_0(y))+\bar\sigma]\leq0.$$
On the other hand, from (\ref{a3.15}) we have, for $0\leq y\leq x_{\bar\varsigma/2}$,
\begin{eqnarray*}
\underline u_t-\underline u_{xx}+\beta \underline u_x-f(\underline u)
           &\leq&Ce^{-\bar\sigma t}\left(\bar\sigma+\max_{0\leq s\leq 1}f'(s)\right)-\kappa C\bar\sigma e^{-\bar\sigma t}(1-Ce^{-\bar\sigma t})V_{\bar\varsigma/2} \\[1mm]
           &=& Ce^{-\bar\sigma t}(1-Ce^{-\bar\sigma t})\left(\frac{\bar\sigma+\max_{0\leq s\leq 1}f'(s)}{1-Ce^{-\bar\sigma t}}-\kappa \bar\sigma V_{\bar\varsigma/2}\right).
\end{eqnarray*}
Hence we can choose  sufficiently large $\kappa>0$ such that $\underline u_t-\underline u_{xx}+\beta \underline u_x-f(\underline u)\leq0$ for $t\geq T$ and $0\leq y\leq x_{\bar\varsigma/2}$.

In short, we derive
\begin{eqnarray*}
\left\{\begin{array}{ll}
\underline u_t-\underline u_{xx}+\beta \underline u_x-f(\underline u)\leq0, & t\geq T,\ \underline l(t)\leq x\leq \underline r(t),\\[1mm]
\underline u(t,\underline l(t))\leq u(t,\underline l(t)),\ \underline u(t,\underline r(t))\leq u(t,\underline r(t)), \ \ &t\geq T,\\[1mm]
\underline u(T,x)\leq  u(T,x), & \underline l(T)\leq x\leq \underline r(T),\\[1mm]
\underline l(t)>0, &t\geq T.
\end{array}\right.
\end{eqnarray*}
The following lemma can be gained by the comparison principle.

\begin{lemma}\label{le7}{\hspace{-1.8mm}\bf .} If $a>0$, then $u(t,x)\geq \underline u(t,x)\ \;{\rm for}\ t\geq T$ and $\underline l(t)\leq x\leq \underline r(t)=c_rt$, where $\underline u(t,x)$ and $\underline l(t)$ are defined by $(\ref{3.15})$ and $\kappa>0$ is suitable large.
\end{lemma}

One knows from Theorem \ref{th2.2} that if spreading happens, then $\lim_{t\to\infty}\|u(t,x)-1\|_{L^\infty_{loc}([0,\infty))}$$=0$ for the case $a=0$, whereas  $\lim_{t\to\infty}\|u(t,x)-\tilde v\|_{L^\infty_{loc}([0,\infty))}=0$  for the case $a>0$, here $\tilde v(x)$ is the unique solution of $(\ref{1.4})$ defined in $[0,\infty)$ and satisfies $\tilde v'(x)>0$, $\tilde v(\infty)=1$. Now we present the uniform convergence for solution $u(t,x)$ and  the sharp estimate for spreading speed of expanding front $h(t)$. To this aim, we first state a proposition concerning the zero set.

Let $\mathcal{Z}_I[u(\cdot)]$ be the number of zeros of a continuous function $u(x)$ defined on $I\subset\mathbb{R}$. The following proposition is an easy consequence of the proofs of \cite[Theorems C and D]{An}.

\begin{proposition}\label{zero2}{\hspace{-1.8mm}\bf .} Let $k>0$ be a constant. Assume that $u(t,x)$ is a bounded classical solution of
\begin{eqnarray*}
\left\{\begin{array}{ll}
u_t=a(t,x)u_{xx}+b(t,x)u_x+c(t,x)u, \ \ &0< t< T,\ 0< x< k,\\[1mm]
u(t,0)=l_0(t),u(t,k)=l_k(t),&0< t< T,
\end{array}
\right.
\end{eqnarray*}
where $l_0(t),l_k(t)\in C^1([0,T])$, and each function is either identically zero or never zero for $t\in[0,T]$. Suppose that $a,1/a,a_t,a_x,a_{xx},b,b_t,b_x,c\in L^\infty$, and $u(0,\cdot)\not\equiv0$ when $l_0=l_k\equiv0$. Then $\mathcal{Z}_{[0,k]}[u(t,\cdot)]<\infty$ for every $t\in(0,T]$ and $\mathcal{Z}_{[0,k]}[u(t,\cdot)]$ is nonincreasing in $t$ for $t\in(0,T]$. Moreover, if for some $t_0\in(0,T]$, $u(t_0,\cdot)$ has a degenerate zero point $x_0\in[0,k]$, then $\mathcal{Z}_{[0,k]}[u(t_1,\cdot)]>\mathcal{Z}_{[0,k]}[u(t_2,\cdot)]$ for any $t_1,t_2$ with $0<t_1<t_0<t_2<T$.
\end{proposition}

The following theorem is the main results of this section.

\begin{theorem}\label{uni}{\hspace{-1.8mm}\bf .} Let $\tilde c_\beta$ and $\tilde q_\beta(z)$ be obtained by Proposition $\ref{p3.1}$.
If $(u,h)$ is a solution of $(P)$ for which spreading happens, then  there exists $H\in\mathbb{R}$ such that
  \bes
  &\ds\lim_{t\to\infty}(h(t)-\tilde c_\beta t-H)=0, \ \ \ \lim_{t\to\infty}h'(t)=\tilde c_\beta;&\lbl{3.16}\\[1mm]
  &\ds\lim_{t\to\infty}\|u(t,x)-\tilde v(x)\tilde q_\beta(\tilde c_\beta t+H-x)\|_{L^\infty([0,h(t)])}=0 \ \ \ {\rm if } \ a>0;&
 \lbl{3.19}\\[1mm]
  &\ds\lim_{t\to\infty}\|u(t,x)- \tilde q_\beta(\tilde c_\beta t+H-x)\|_{L^\infty([0,h(t)])}=0 \ \ \ {\rm if } \ a=0.&
 \lbl{3.20a}\ees
Here we have employed the convention that $\tilde q_\beta$ is extended to be zero outside its support.
\end{theorem}

\noindent{\bf Proof.} The argument will be divided into the following four steps. In the first step we shall prove (\ref{3.16}). To show the conclusions (\ref{3.19}) and (\ref{3.20a}), the locally uniform convergence of $u$ near $h(t)$  is discussed in the second step. The last two steps are devoted to (\ref{3.19}) and (\ref{3.20a}), respectively.

{\it Step 1: Proof of (\ref{3.16})}. Define
 \bess
 &\tilde g(t)=-\tilde c_\beta t,\ \ \tilde s(t)=h(t)-\tilde c_\beta t\ \ \ {\rm for }\ \ t\geq0,&   \\
 &\zeta(t,y)=u(t,y+\tilde c_\beta t)\ \ \ {\rm for}\ \ t\geq0,\ \ y\in[\tilde g(t),\tilde s(t)].&
  \eess
Then $(\zeta,\tilde g, \tilde s)$ solves the following problem:
\begin{equation}\label{a3.17}
\left\{
\begin{array}{ll}
\zeta_t-\zeta_{yy}-(\tilde c_\beta-\beta)\zeta_y=f(\zeta),&t>0,\ \tilde g(t)<y<\tilde s(t),\\[1mm]
B[\zeta](t,\tilde g(t))=0,\ \tilde g'(t)=-\tilde c_\beta,&t>0,\\[1mm]
\zeta(t,\tilde s(t))=0,\ \tilde s'(t)=-\mu\zeta_y(t,\tilde s(t))-\tilde c_\beta,\ \ & t>0,\\[1mm]
\tilde g(0)=0,\ \ \tilde s(0)=h_0,\ \ \zeta(0,y)=u_0(y),&0\leq y\leq h_0.
\end{array}\right.
\end{equation}

For every $d\in\mathbb{R}$, the function $v(y)=\tilde q_\beta(d-y)$ is a stationary solution of the equation in (\ref{a3.17}):
\begin{eqnarray*}
\left\{
\begin{array}{l}
v''+(\tilde c_\beta-\beta)v'+f(v)=0,\ \ -\infty<y<d,\\[1mm]
v(-\infty)=1,\ \ v(d)=0,\ \ v'(d)=-\tilde c_\beta/\mu.
\end{array}\right.
\end{eqnarray*}
By the same argument as in \cite[Lemma 3.7]{KM} to study the number of zeroes of $\zeta(t,\cdot)-v(\cdot)$ one can derive that $\tilde s(t)-d$ changes its sign at most finite many times. On the other hand, Lemma \ref{le5} gives  $\underline h(t)\leq h(t)\leq \bar h(t)$ for sufficiently large $t$, which implies the boundedness of $h(t)-\tilde c_\beta t$. Therefore, there are  a sequence $\{t_n\}\in \mathbb{R}$ satisfying $t_n\to\infty$ as $n\to\infty$ and $H\in\mathbb{R}$ such that $\tilde s(t_n)\to H$ as $n\to \infty$. Assume that there exist another sequence $\{\tilde t_n\}\in\mathbb{R}$ satisfying $\tilde t_n\to\infty$ as $n\to\infty$ and $\tilde H\neq H$ such that $\tilde s(\tilde t_n)\to \tilde H$ as $n\to \infty$. So $\tilde s(t)-d$ changes its sign infinite many times for $\min\{H,\tilde H\}<d<\max\{H,\tilde H\}$. This is a contradiction, and thus $\lim_{t\to\infty}\tilde s(t)=H$.

It follows from the estimate (\ref{1.3}) that $\tilde s'(t)$ is H\"older continuous and $\|\tilde s'\|_{C^{\alpha}([n+1,n+3])}\leq C$ for some $C>0$, $\alpha\in(0,1)$ and all $n\ge 1$. By virtue of the limit of $\tilde s(t)$ it is not difficult to verify that $\lim_{t\to\infty}\tilde s'(t)=0$, which implies $\lim_{t\to\infty}h'(t)=\tilde c_\beta$.

{\it Step 2: Locally uniform convergence of $u(t,y+h(t))$}.
Set
$$\nu(t)=-h(t)\ \ {\rm for}\ t\geq0;\ \ \ \xi(t,y)=u(t,y+h(t)) \ \ {\rm for}\ t\geq0,\ \ y\in[\nu(t),0].$$
Then $(\xi(t,x),\nu(t))$ is a solution of
\begin{eqnarray*}
\left\{
\begin{array}{ll}
\xi_t-\xi_{yy}-(h'(t)-\beta)\xi_y=f(\xi),&t>0,\ \nu(t)<y<0,\\[1mm]
B[\xi](t,\nu(t))=0,\ \nu'(t)=\mu\xi_y(t,0),\ \ &t>0,\\[1mm]
\xi(t,0)=0,\ h'(t)=-\mu\xi_y(t,0),& t>0,\\[1mm]
\nu(0)=-h_0,\ \ \xi(0,y)=u_0(y+h_0),&-h_0\leq y\leq 0.
\end{array}\right.
\end{eqnarray*}
Obviously, $\nu_\infty=\ds\lim_{t\to\infty}\nu(t)=-\infty$. By the parabolic $L^p$ theory and imbedding theorem we affirm that for any sequence $\{t_n\}_{n\in \mathbb{N}}$ satisfying $t_n\to\infty$ as $n\to\infty$, there exists a subsequence $\{t_n'\}\subset\{t_n\}$ such that
$$\xi(t+t_n',y)\to v(t,y)\ \ {\rm as }\ n\to\infty {\ \rm locally\ uniformly\ in\ } (t,y)\in\mathbb{R}\times(-\infty,0],$$
and $v(t,y)$ satisfies
\begin{eqnarray*}
\left\{
\begin{array}{ll}
v_t-v_{yy}-(\tilde c_\beta-\beta) v_y=f(v),&t\in\mathbb{R},\ y<0,\\[1mm]
v(t,0)=0,\ v_y(t,0)=-\tilde c_\beta/\mu, \ \ &t\in\mathbb{R}.
\end{array}\right.
\end{eqnarray*}

We assert that $ v(t,y)\equiv \tilde q_\beta(-y)$ for every $t\in\mathbb{R}$. Assume on the contrary that there exists $(t_0,y_0)\in\mathbb{R}\times(-\infty,0)$ such that $v(t_0,y_0)\neq \tilde q_\beta(-y_0)$. Due to the continuity of $v(t,y)$, there exists $0<\varepsilon\ll 1$ such that $v(t_0+t,y_0)\neq \tilde q_\beta(-y_0)$ for all $t\in(0,\varepsilon)$. Applying the zero number result (Proposition \ref{zero2}) to  $\eta(t,y)=v(t_0+t,y)-\tilde q_\beta(-y)$ in $[0,\varepsilon]\times[y_0,0]$, one can easily show that the number of zeroes of $\eta(t,y)$ in $[y_0,0]$ is finite for $t\in(0,\varepsilon)$, and it decreases strictly once it has a degenerate zero point in $[y_0,0]$. It is easy to derive that $\eta(t,0)=\eta_y(t,0)\equiv0$ for $t\in(0,\varepsilon)$, i.e., $y=0$ is a degenerate zero point of $\eta(t,\cdot)$ for each  $t\in(0,\varepsilon)$. We get a contradiction.

As a consequence, $u(t,y+h(t))\to\tilde q_\beta(-y)$ as $t\to\infty$ uniformly on $[-L,0]$ for any given $L>0$.

{\it Step 3: Proof of (\ref{3.19})}.
Fix $c\in(0,\tilde c_\beta)$. We first show that
 \bes
 \lim_{t\to\infty}\|u(t,x)-\tilde v(x)\|_{L^\infty([0,\,c_rt])}=0,
 \lbl{3.20}\ees
 where $c_r$ is defined by (\ref{alr}).
Thanks to the second limit of (\ref{3.16}) and $c<\tilde c_\beta$, we have that $h(t)> c_rt$ for $t\gg 1$. In view of Lemmas \ref{le6} and \ref{le7}, there exists a constant $T\gg 1$ such that, for $t\ge T$ and $x\in[0,c_rt]$,
$$(1-Ce^{-\bar\sigma  t})\tilde v_0(x-\underline l(t))\leq u(t,x)\leq (1+Me^{-\sigma  t})\tilde v(x-\bar l(t)),$$
where we have assumed that $\tilde v_0(z)=0$ for $z\leq0$, and
 \[\bar l(t)=\kappa M(e^{-\sigma t}- e^{-\sigma T})-x^*, \ \ \underline l(t)=c_lT-\kappa C(e^{-\bar\sigma t}- e^{-\bar\sigma T}).\]
Let $\tilde v(z)=0$ for $z\leq0$. Noting $f'(1)<0$, a standard argument generates that there exist positive constants $K$ and $\rho$ such that
 \bes\label{a3.18}  1-\tilde v(z)\leq K e^{-\rho z}, \ \ 1-\tilde v_0(z)\leq K e^{-\rho z},\ \ \forall \ z\in\mathbb{R}. \ees
Applying (\ref{a3.18}) and the boundedness of $\bar l(t)$ and $\underline l(t)$,  we deduce that, for $t\ge T$ and $0\le x\le c_rt$,
\begin{eqnarray}
1-u(t,x)&\leq& 1-(1-Ce^{-\bar\sigma  t})\tilde v_0(x-\underline l(t))\nonumber\\[.1mm]
        &=&1-\tilde v_0(x-\underline l(t))+Ce^{-\bar\sigma  t}\tilde v_0(x-\underline l(t)\nonumber\\[.1mm]
        &\leq&K e^{-\rho (x-\underline l(t))}+Ce^{-\bar\sigma  t}\nonumber\\[.1mm]
        &\leq& K'(e^{-\rho x}+e^{-\bar\sigma t}),\label{a3.19}
\end{eqnarray}
and
 \bes
u(t,x)-1&\leq& (1+Me^{-\sigma  t})\tilde v (x-\bar l(t))-1\nonumber\\[.1mm]
        &\leq&1-\tilde v (x-\bar l(t))+Me^{-\sigma  t}\tilde v (x-\bar l(t)\nonumber\\[.1mm]
        &\leq&K e^{-\rho (x-\bar l(t))}+Me^{-\sigma t}\nonumber\\[.1mm]
        &\leq& K'(e^{-\rho x}+e^{-\bar\sigma t})\label{a3.20}
\ees
on account of $\bar\sigma<\sigma$. Combining (\ref{a3.18}), (\ref{a3.19}) and (\ref{a3.20}), we see that, for any small $\varepsilon>0$, there exist large positive constants $C$ and $\hat T$ so that
 \bes
 \sup_{x\in[C,\,c_rt]}|u(t,x)-\tilde v(x)|<\varepsilon, \ \ \forall \ t>\hat T.
 \lbl{3.24}\ees

On the other hand, by the locally uniform convergence of $u(t,x)$ to $\tilde v(x)$ on $[0,\infty)$ (cf. Theorem \ref{th2.2}(i)), we obtain
  $$\sup_{x\in[0,C]}|u(t,x)-\tilde v(x)|<\varepsilon\ \ {\rm for\ }t\gg 1.$$
This combined with (\ref{3.24}) allows us to derive (\ref{3.20}).

Next it will be proved that
\begin{equation}\label{a3.21}
{\displaystyle\lim_{t\to\infty}\|u(t,x)-\tilde q_\beta(\tilde c_\beta t+H-x)\|_{L^\infty([c_lt, h(t)])}=0},\end{equation}
where $c_l$ is defined by (\ref{alr}).
By Lemma \ref{le5}, there exists $T^*>0$ such that, for $t\ge T^*$ and $x\in[c_lt,h(t)]$,
  \begin{equation}\label{a3.22}
  (1-K_2e^{-\bar\sigma  t})\tilde q_\beta(\underline h(t)-x)\leq u(t,x)\leq (1+K_1e^{-\sigma  t})\tilde q_\beta(\bar h(t)-x),\end{equation}
where $\bar h(t)$ and $\underline h(t)$ are given by (\ref{3.6}) and (\ref{3.7}), respectively. Employing the same argument to (\ref{a3.18}) yields that, for some $K$, $\rho>0$,
  \bes\label{a3.23}
1-\tilde q_\beta(z)\leq K e^{-\rho z}, \ \ \ \forall\ z\in\mathbb{R}.\ees
Utilizing the first conclusion of (\ref{3.16}) and the expressions of $\bar h(t)$ and $\underline h(t)$, it is easy to see that $\underline h(t)-h(t)$ and $\bar h(t)-h(t)$ are bounded. Thus,
by (\ref{a3.22}) and (\ref{a3.23}), there exists a positive constant $K'$ such that
   $$|1-u(t,y+h(t))|\leq K'(e^{\rho y}+e^{-\bar\sigma t}), \ \ \ \forall \ t\ge T^*,\ \, y\in[c_lt-h(t),0]. $$
Thus, for any small $\varepsilon>0$, there exist large positive constants $C$ and $T'$ so that
  $$\sup_{y\in[c_lt-h(t),\,-C]}|u(t,y+h(t))-\tilde q_\beta(-y)|<\varepsilon, \ \ \ \forall \ t>T'.$$
On the other hand, it follows from Step 2 that
$$\sup_{y\in[-C,\,0]}|u(t,y+h(t))-\tilde q_\beta(-y)|<\varepsilon \, \ \ {\rm for\,\ }t\gg 1.$$
Consequently, we obtain
  $$\sup_{y\in[c_lt-h(t),\,0]}|u(t,y+h(t))-\tilde q_\beta(-y)|<\varepsilon\, \ \ {\rm for\,\ }t\gg 1,$$
which implies
\begin{equation}\label{aks}\sup_{x\in[c_lt,\,h(t)]}|u(t,x)-\tilde q_\beta(h(t)-x)|<\varepsilon\, \ \ {\rm for\,\ }t\gg 1.\end{equation}
Combining (\ref{aks}) and the first limit of (\ref{3.16}) enables us to deduce (\ref{a3.21}).

The conclusion (\ref{3.19}) can be achieved from (\ref{3.20}) and (\ref{a3.21}). In fact, for any small $\varepsilon>0$ and $x\in[0,c_rt]$,
$$1-\tilde q_\beta(\tilde c_\beta t+H-x)<\varepsilon \ \  {\rm for\,\ }t\gg 1 $$
because $\tilde c_\beta t+H-x\geq c_lt\to\infty$ as $t\to\infty$. Combining this with (\ref{3.20}), we can get,
for any $x\in[0,c_rt]$,
\begin{eqnarray*}&&|u(t,x)-\tilde v(x)\tilde q_\beta(\tilde c_\beta t+H-x)|\leq |u(t,x)-\tilde v(x)|+\tilde v(x)|1-\tilde q_\beta(\tilde c_\beta t+H-x)|\\
&&\leq |u(t,x)-\tilde v(x)|+1-\tilde q_\beta(\tilde c_\beta t+H-x)\leq \varepsilon+\varepsilon=2\varepsilon \ \ \  {\rm for\,\ }t\gg 1.
\end{eqnarray*}
Analogously, for the sufficiently large $t$ and any $x\in[c_lt,h(t)]$,
\[|u(t,x)-\tilde v(x)\tilde q_\beta(\tilde c_\beta t+H-x)|\leq|u(t,x)-\tilde q_\beta(\tilde c_\beta t+H-x)|+1-\tilde v(x)\leq \varepsilon.\]

{\it Step 4: Proof of (\ref{3.20a})}. In this case we have $B[u](0)=u_x(t,0)=0$. Take $\underline g(t)=0$ and $\kappa\gg 1$ in (\ref{3.7}). It is easy to verify that
  \begin{eqnarray*}
\left\{\begin{array}{ll}
\underline u_t-\underline u_{xx}+\beta \underline u_x-f(\underline u)\leq0, & t\geq T,\ 0\leq x\leq \underline h(t),\\[1mm]
\underline u_x(t,0)\leq 0,\ \underline u(t,\underline h(t))=0,\ \ &t\geq T, \\[1mm]
\underline h'(t)\leq -\mu \underline u_x(t,\underline h(t)),&t\geq T,\\[1mm]
\underline u(T,x)\leq  u(T,x), & 0\leq x\leq \underline h(T)<h(T).
\end{array}\right.
  \end{eqnarray*}
By the comparison principle, we derive that $u(t,x)\geq\underline u(t,x)$ for $t\geq T$ and $0\leq x\leq \underline h(t)$. On the other hand,  $(\bar u,\bar g,\bar h)$ in (\ref{3.6})  is also an upper solution of $(P)$ for this case.  Thus (\ref{a3.22}) still holds for $t\geq T^*$ and $0\leq x\leq h(t)$ with $T^*$ sufficiently large. Applying the same way as (\ref{a3.21}) we can reach the desired results.
\ \ \ \fbox{}

\begin{remark}{\hspace{-1.8mm}\bf .}
From Theorem \ref{uni} and Proposition \ref{p3.1}, it is easy to see that the spreading speed $\tilde c_\beta$ is  getting slower as $\beta$ becomes small in $(-c_0,c_0)$, and approaches $0$ as $\beta$ tends to $-c_0$.
\end{remark}

\section{Long time behavior of solutions for either $\beta\geq c_0$ and $a\geq bc_0/2$ or $\beta\leq- c_0$}\label{se4}
\setcounter{equation}{0}

In this section we only present a brief description concerning the long time behavior of solutions $u(t,x)$ to problem $(P)$ in cases either $\beta\geq c_0$ and $a\geq bc_0/2$ or $\beta\leq- c_0$.
Firstly, a locally uniform convergence conclusion is provided for the cases.

\begin{theorem}\label{t4.1}{\hspace{-1.8mm}\bf .} If either $\beta\geq c_0$ and $a\geq bc_0/2$, or $\beta\leq- c_0$, then the solution $u(t,\cdot)$ of problem $(P)$ converges to $0$ locally uniformly in $[0,h_\infty)$ as $t\to\ty$ regardless of $h_\infty<\infty$ or $h_\infty=\infty$.
\end{theorem}

\noindent{\bf Proof.} When $\beta<-c_0$, we first consider the following problem
\begin{eqnarray}\label{a4.1}
\left\{
\begin{array}{l}
q''(z)-cq'(z)+f(q)=0,\ \ z\in\mathbb{R},\\[1mm]
q(-\infty)=1,\ q(\infty)=0,\\[1mm]
 q'(z)<0,\ \ z\in\mathbb{R}.
\end{array}\right.
\end{eqnarray}
It is well known (see, for example, \cite{AW,KPP}) that problem (\ref{a4.1}) admits a solution $q(z;c)$ if and only if $c\leq -c_0$. Set $q^-(z)=q(z;-c_0)$, then $u(t,x)=q^-(x-(\beta+c_0)t)$ is a traveling wave of $u_t-u_{xx}+\beta u_x=f(u)$, and  travels leftward  if and only if  $\beta<-c_0$.

Let $k=2\max\{1,\|u_0\|_{L^\infty([0,h_0])}\}$, and
define a function $f_k(s)\in C^2([0,\infty))$ satisfying
\bess
f_k(s)
\left\{\begin{array}{ll}=f'(0)s,\ \ & 0\leq s\leq 1,\\[.1mm]
>0,& 0< s < k,\\[.1mm]
<0,& s>k,
 \end{array}\right.\eess
and
 \[f'_k(k)<0,\ \ f(s)\leq f_k(s)\leq f'_k(0)s \ \ \ {\rm for }\, \ s\geq0.\]
Denote by $q^-_k(z)$ the unique solution of (\ref{a4.1}) with
$c=-c_0$, with $f$ and $q(-\infty)=1$ replaced by $f_k$ and $q(-\infty)=k$, respectively. Notice that
  \[B[q^-_k]\bbb(-(\beta+c_0)t-x_0\bbb)\geq0 \ \ \ {\rm for}\, \ t>0,\]
and
  \[u_0(x)\leq q^-_k(x-x_0)\ \ \ {\rm for}\, \ x\in[0,h_0]\]
provided that $x_0>0$ is sufficiently large. The comparison principle yields
\begin{equation}\label{a4.2}
u(t,x)\leq q^-_k\bbb(x-(\beta+c_0)t-x_0\bbb)\ \ \  {\rm for }\ \ t>0, \ 0\leq x\leq h(t).
\end{equation}
As a consequence, $u(t,x)$ converges to $0$ locally uniformly in $x\in[0, h_\infty)$ as $t\to\ty$ since $q^-_k(x-(\beta+c_0)t-x_0)$
 is a leftward traveling wave with the positive speed $-(\beta+c_0)$ and $q^-_k(z)\to0$ as $z\to\infty$.

When $\beta=-c_0$, using a similar argument as in \cite{GLZ}, one can show that there exists a positive constant $M$ such that
\begin{equation}\label{a4.3}
u(t,x)\leq Me^{-\frac{5c_0}{12}(x+p(t))}\ \  {\rm for}\ t>0,\ \max\bbb\{0,h_0-p(t)\bbb\}\leq x\leq\min\bbb\{h(t),\sqrt{t+1}-p(t)\bbb\},
\end{equation}
where $p(t)=\frac3{c_0}\ln(1+\frac{t}{t_0})$ for $t>0$. Then for any $0<L<h_\infty$, when $t$ is sufficiently large we have
   $$L<\sqrt{t+1}-p(t),$$
Thus, by virtue of (\ref{a4.3}), for any $x\in[0,L]$,
$$u(t,x)\leq Me^{-\frac{5c_0}{12}(L+p(t))}\to0\ \ \;{\rm as}\ t\to\infty.$$

If $\beta\geq c_0$, similar to the above discussion, the problem
\begin{eqnarray}\label{a4.4}
\left\{
\begin{array}{l}
q''(z)-cq'(z)+f(q)=0,\ \ z\in\mathbb{R},\\[1mm]
q(-\infty)=0,\ q(\infty)=1,\\[1mm]
 q'(z)>0,\ \ z\in\mathbb{R}
\end{array}\right.
\end{eqnarray}
 admits a solution $q(z;c)$ if and only if $c\geq c_0$.  Denote $q^+(z)=q(z;c_0)$, then the differential equation of $(P)$ admits a traveling wave $u(t,x)=q^+(x-(\beta-c_0)t)$, which travels rightward if and only if $\beta> c_0$. Denote by $q^+_k(z)$ the unique solution of (\ref{a4.4}) with
$c=c_0$, with $f$ and $q(\infty)=1$ replaced by $f_k$ and $q(\infty)=k$, respectively. It should be emphasized that $B[q^+_k]\big(-(\beta-c_0)t+x_0\big)\geq0$ for all $t>0$ when $a\geq bc_0/2$.  We can  also deduce the estimates corresponding to (\ref{a4.2}) and (\ref{a4.3}):
\begin{equation}\label{b4.5}u(t,x)\leq q^+_k(x-(\beta-c_0)t+x_0) \ \ \ {\rm for }\ \ t>0, \ \ 0\leq x\leq h(t)
\end{equation}
provided that $x_0>0$ is sufficiently large, and
$$
u(t,x)\leq Me^{-\frac{5c_0}{12}(p(t)-x)} \ \ \ {\rm for}\ t>0,\ \max\{0,p(t)-\sqrt{t+1}\}\leq x\leq\min\{h(t),p(t)-h_0\}.
$$
Similarly to the above, the required conclusion for $\beta\geq c_0$ can be obtained from these two estimates. The proof is complete.
\ \ \ \fbox{}

\begin{remark}{\hspace{-1.8mm}\bf .} The above argument concerning the case $\beta\leq-c_0$ enables us to conclude that $u(t,x)$ uniformly converges to $0$ on $[0,h_\infty)$ as $t\to\infty$. That is to say, the species will extinct eventually for $\beta\leq-c_0$ whether $h_\infty<\infty$ or $h_\infty=\infty$.
\end{remark}

Now we present the boundedness of expanding front $h(t)$ when $\beta<-c_0$.

\begin{theorem}\label{t4.2}{\hspace{-1.8mm}\bf .} Assume $\beta<-c_0$ and $(u,h)$ is a solution of problem $(P)$. Then $h_\infty<\infty$.
\end{theorem}

\noindent{\bf Proof.} Let $x_0>0$ be such that (\ref{a4.2}) holds. Note that
\begin{equation}\label{b4.6}q^-_k(y)\sim C ye^{-\frac{c_0}2y}\ \ \ {\rm as}\ \ y\to\infty\end{equation}
for some positive constant $C$ (see \cite{AW,HNRR}), and $\beta+c_0<0$, we can find $T,\,M>0$ such that, for $t\geq T$ and $h_0\leq x\leq h(t)$,
\bess
\begin{array}{ll}
q^-_k\big(x-(\beta+c_0)t-x_0\big)&\leq q^-_k\big(h_0-(\beta+c_0)t-x_0\bbb)\\[1mm]
&\leq C\bbb(h_0-(\beta+c_0)t-x_0\bbb)e^{-\frac{c_0}2\bbb(h_0-(\beta+c_0)t-x_0\bbb)}\\[1mm]
&\leq M e^{\frac{c_0}4(\beta+c_0)t}.
\end{array}
\eess
Thus we get, by  (\ref{a4.2}),
\begin{eqnarray*}
u(t,x)\leq q^-_k(x-(\beta+c_0)t-x_0)\leq M t e^{\frac{c_0}4(\beta+c_0)t}
\end{eqnarray*}
for $t\geq T$ and $h_0\leq x\leq h(t)$. Set $\delta=\min\{1,-\frac{c_0}4(\beta+c_0)\}$, $\varepsilon=Me^{\frac{c_0(\beta+c_0)T+\beta\pi}{4}}$,
  $$g(t)=h(T)+\frac\pi2+\frac{\mu\varepsilon}\delta(1-e^{-\delta t})\ \ \ {\rm for} \ t\geq0,$$
and
 $$w(t,x)=\varepsilon e^{-\delta t}e^{-\frac\beta2(g(t)-x)}\sin(g(t)-x)\ \ \ {\rm for}\ \ t\geq0, \ g(t)-\frac\pi2\leq x\leq g(t).$$
By a series of calculations one can verify that, for $t\geq0$ and $g(t)-\frac\pi2\leq x\leq g(t)$,
\begin{eqnarray*}
\begin{array}{lll}
\displaystyle w_t- w_{xx}+\beta w_x-f( w)&= & \displaystyle w\left(-\delta-\frac\beta2g'(t)+\frac{\beta^2}4+1\right)-f(w)\\[1mm]
 &&+ g'(t)\varepsilon e^{-\delta t}e^{-\frac\beta2(g(t)-x)}
 \cos(g(t)-x)\\[1mm]
&\geq&\displaystyle w\left(-\delta-\frac\beta2\mu\varepsilon+\frac{\beta^2}4+1-f'(0)\right)\geq0.
\end{array}
\end{eqnarray*}
On the other hand, it is easy to see that $-\mu w_x(t,g(t))={\mu\varepsilon}e^{-\delta t}=g'(t)\ \ \ {\rm for} \ \, t\geq0,$ and
  $$w(t,g(t)-\frac\pi2)=\varepsilon e^{-\delta t}e^{-\frac{\beta\pi}4}\geq Me^{\frac{c_0(\beta+c_0)}{4}(t+T)}\geq u(t+T,x) \ \ {\rm for} \ t\geq0, \  h_0\leq x\leq h(t).$$
Therefore, if $h(t)$ and $g(t)-\frac\pi2$ do not intersect for all $t>0$, then $h(t)<g(t)-\frac\pi2\leq h(T)+\frac{\mu\varepsilon}\delta$; otherwise, $u(t+T,x)$ and $w(t,x)$ have a common domain, and then the comparison principle generates
$$h(t)\leq g(t)\leq h(T)+\frac\pi2+\frac{\mu\varepsilon}\delta<\infty \ \ {\rm for}\ t>0.$$
This illustrates $h_\infty<\infty$ for $\beta<-c_0$.
\ \ \ \fbox{}
\vspace{.3cm}

Next we further consider the long time behavior of solutions for the case $\beta\geq c_0$ and $a\geq bc_0/2$, which is essentially identical to that of  $(\ref{1.1})$ with $\beta\geq c_0$ obtained in \cite{GLZ} because the solution locally converges to zero as the time $t$ tends to infinity by Theorem \ref{t4.1}. For the completeness and for the convenience to the reader, we present a summary and brief discussion. Before stating the conclusions for $(P)$ with advection $\beta\geq c_0$, several concepts (see \cite{GLZ}) are listed as follows:
\begin{itemize}
\item {\it virtual spreading }: $h_\infty=\infty$, $\ds\lim_{t\to\infty}u(t,x)=0$ locally uniformly in $[0,\infty)$, and
 $$\lim_{t\to\infty}u(t,x+ct)=1\ \,{\rm locally\ uniformly\ in}\ \mathbb{R},\ {\rm for\ some} \ c>0;$$
\item {\it vanishing }: $h_\infty<\infty$ and
\begin{equation}\label{bbb}
\lim_{t\to\infty}\max_{0\leq x\leq h(t)}u(t,x)=0;
\end{equation}
\item {\it virtual vanishing }: $h_\infty=\infty$, and (\ref{bbb}) holds.
\end{itemize}

Firstly, the conclusions for the problem with advection $\beta\geq\beta^*$ is provided. Here and in what follows, $\beta^*$ (see Remark \ref{k3.1}) is the unique root of the equation $\beta-c_0=\tilde c_\beta$.

\begin{theorem}{\hspace{-1.8mm}\bf .}\label{t4.3} Assume that $\beta\geq\beta^*$, $a\geq bc_0/2$, and $(u,h)$ is a time-global solution of $(P)$. Then vanishing happens {\rm($h_\infty<\infty$)}.
\end{theorem}
\noindent{\bf Proof.} The argument is almost identical to that of Proposition 4.7 in \cite{GLZ}, so we only give a sketch here for the sake of convenience.

We first discuss the case $\beta>\beta^*$ and $a\geq bc_0/2$. In this case we have $\tilde c_\beta<\beta-c_0$. Denote $\vartheta=\beta-c_0-\tilde c_\beta>0$.  Note the results of Lemma \ref{le5} still hold for $\beta\geq c_0$. Similarly to the estimate (\ref{b4.6}), one can know
\begin{equation}\label{b4.7}q^+_k(y)\sim -C ye^{\frac{c_0}2y}\ \ \ {\rm as}\ \ y\to-\infty\end{equation}
for some positive $C$. Therefore, in view of (\ref{b4.5}), (\ref{b4.7}) and Lemma \ref{le5}, there exist $T_1>0$, $C_1>0$ such that, for $t\geq T_1$ and $x\in[0,h(t))$,
\begin{eqnarray*}
\begin{array}{ll}
u(t,x)&\leq q^+_k\big(x-(\beta-c_0)t+x_0\big)\leq q^+_k\big(h(t)-(\beta-c_0)t+x_0\bbb)\\[1mm]
&\leq  q^+_k\big(-\vartheta t+M_1+x_0\bbb)\leq -2C(-\vartheta t+M_1+x_0) e^{\frac{c_0}2(-\vartheta t+M_1+x_0)}\\[1mm]
&\leq C_1 e^{-\frac{c_0\vartheta}4t}.
\end{array}
\end{eqnarray*}
Let $\delta=\frac12\min\{1,c_0\vartheta\}$ and select $T_2>T_1$ so that
$$\varepsilon=C_1 e^{\frac{\beta\pi-c_0\vartheta T_2}4}<\frac2{\beta\mu}.$$
Define
$$g(t)=h(T_2)+\frac\pi 2+\frac{\mu\varepsilon}{\delta}(1-e^{-\delta t})\ \,{\rm for}\ t\geq0$$
and $$w(t,x)=\varepsilon e^{-\delta t} e^{\frac\beta 2(x-g(t))}\cos\left(x-g(t)+\frac\pi2\right)\ \,{\rm for}\ t\geq0,\ g(t)-\frac\pi2\leq x\leq g(t).$$
A straightforward calculation as in the proof of Theorem \ref{t4.2} concludes that $(w(t,x),g(t)-\frac\pi2,g(t))$ is an upper solution and
$$h(t+T_2)\leq g(t)\leq h(T_2)+\frac\pi2+\frac{\mu\varepsilon}\delta<\infty.$$

Now we consider the case $\beta=\beta^*$ and $a\geq bc_0/2$. Suppose on the contrary that $h_\infty=\infty$, then by Theorem \ref{t4.1},  $u(t,x)\to 0$  ($t\to\infty$) on $x\in[0,L]$ for every positive constant $L$. Thus we can choose $T_3>0$ and $L'>0$ such that $u(T_3,x)$ and $\tilde q_{\beta^*}(\tilde c_{\beta^*}T_3-L'-x)$ intersect only one point at the right small neighbourhood of $x=0$, where $\tilde c_{\beta^*}$ and $\tilde q_{\beta^*}$ are established in Remark \ref{k3.1}. Notice that $\tilde q_{\beta^*}(\tilde c_{\beta^*}(t+T_3)-L'-x)$ is the rightward traveling semi-wave with end point at $\tilde c_{\beta^*}T_3-L'$.
 Next, utilizing a completely similar argument as in \cite[Prosition 4.7]{GLZ}, we can display that, for some large $T_4>T_3$,
$$u(t,x)<\tilde q_{\beta^*}(\tilde c_{\beta^*}t-L'-x)\ \,{\rm for}\ t\geq T_4,\ 0\leq x\leq h(t),$$
and deduce a contradiction with the indirect assumption.
\ \ \ \fbox{}

For problem $(P)$ with advection $c_0\leq \beta<\beta^*$, we have the following conclusions, which are nearly identical to Theorems 2.2 and 2.3 of \cite{GLZ}.

\begin{theorem}{\hspace{-1.8mm}\bf .}\label{t4.4} Assume $\beta=c_0$, $a\geq b\beta/2$, and $(u,h)$ is a time-global solution of $(P)$ with $u_0=\lambda \psi$ and $\psi\in\mathscr{X}(h_0)$. Then there exist $\lambda_*$, $\lambda^*\in(0,\infty]$ with $\lambda_*\leq\lambda^*$ such that

{\rm (i)} if $\lambda>\lambda^*$, then virtual spreading happens for any $c\in(0,\tilde c_\beta)$, where $\tilde c_\beta$ is given in Remark \ref{k3.1};

{\rm (ii)} if $0<\lambda<\lambda_*$, then vanishing happens;

{\rm (iii)} if $\lambda_*\leq\lambda\leq\lambda^*$, then virtual vanishing happens.
\end{theorem}

\begin{theorem}{\hspace{-1.8mm}\bf .}\label{t4.5} Assume that $c_0<\beta<\beta^*$, $b=0$, and $(u,h)$ is a time-global solution of $(P)$ with $u_0=\lambda \psi$ and $\psi\in\mathscr{X}(h_0)$. Then there exists $\lambda_*\in(0,\infty]$ dependent on $h_0, \psi$ such that

{\rm (i)} if $\lambda>\lambda_*$, then virtual spreading happens for any $ c\in(\beta-c_0,\tilde c_\beta)$, where $\tilde c_\beta$ is given in Remark \ref{k3.1};

{\rm (ii)} if $0<\lambda<\lambda_*$, then vanishing happens;

{\rm (iii)} if $\lambda=\lambda_*$, then $\ds\lim_{t\to\infty}h'(t)=\beta-c_0$, and
$$\lim_{t\to\infty}\|u(t,x)-V^*(x-h(t))\|_{L^\infty([0,h(t)])}=0,$$
where $V^*(z)$ is the unique solution of
\begin{eqnarray}\label{b4.9}
\left\{\begin{array}{l}
V''(z)-c_0V'(z)+f(V)=0, \ \ V(z)>0,\ \ z\in(-\infty,0),\\[1mm]
V(0)=0,\ \ V(-\infty)=0, \ \ -\mu V'(0)=\beta-c_0.\end{array}\right.
\end{eqnarray}
\end{theorem}

\begin{remark}In \cite{GLZ}, $V^*$ is called a tadpole-like shape: it has a ``big head''  on the right side and an infinite lone ``tail'' on the left side. Moreover, $V^*(x-(\beta-c_0)t)$ is called a tadpole-like traveling wave with speed $\beta-c_0$, which exists if and only if $\beta\in(c_0,\beta^*)$(see \cite[Lemma 3.5]{GLZ}).
\end{remark}

In order to clarify the above two theorems, we first show that the solution $u(t,x)$ of $(P)$ with $\beta\geq c_0$ converges uniformly to zero on $[0,h_\infty)$ when the initial data $u_0$ is sufficiently small.
\begin{lemma} Assume that $\beta\geq c_0$, $a\geq b\beta/2$, $(u,h)$ is the solution of problem $(P)$. If $\|u_0\|_{L^\infty([0,h_0))}$ is sufficiently small, then vanishing happens.
\end{lemma}
\noindent{\bf Proof.}
We can select $\delta>0$ small enough such that
$$\frac{\pi^2}{h_0^2(1+\delta)^2}\geq 4\delta+\beta h_0\delta^2.$$
Let $g(t)=h_0(1+\delta-\frac\delta2 e^{-\delta t})$, $\varepsilon=\frac{h_0^2\delta^2}{\pi\mu}(1+\frac\delta2)$ and
$$w(t,x)=\varepsilon e^{-\delta t}e^{\frac\beta2(x-g(t))}\cos\frac{\pi x}{2g(t)}\ \,{\rm for}\ t>0,\ 0\leq x\leq g(t).$$
Similarly to the computation of Lemma \ref{lb.7}, one can achieve
 \begin{eqnarray*}
\left\{\begin{array}{ll}
w_t-w_{xx}+\beta w_x-f(w)\geq0, & t> 0,\ 0\leq x\leq g(t),\\[1mm]
B[w](t,0)\geq 0, \ \ w(t,g(t))=0, \ \ &t> 0,\\[1mm]
g'(t)\geq -\mu w_x(t,g(t)),&t>0.
\end{array}\right.
  \end{eqnarray*}
Take $\|u_0\|_{L^\infty([0,h_0))}$ small enough such that
$$u_0(x)\leq w(0,x)\ \,{\rm for}\ x\in[0,h_0].$$
Therefore, the comparison principle yields $h(t)\leq g(t)\leq h_0(1+\delta)$ for $t>0$ and
\begin{eqnarray*}
\|u(t,\cdot)\|_{L^\infty([0,h(t))}\leq \|w(t,\cdot)\|_{L^\infty([0,g(t))}\leq \varepsilon e^{-\delta t}\to 0\ \,{\rm as}\ t\to\infty.
\end{eqnarray*}
The proof is complete.
\ \ \ \fbox{}

Applying the continuity method and Theorem \ref{t4.1} (also see \cite[Theorem 4.9]{GLZ}), one can derive
\begin{proposition}{\hspace{-1.8mm}\bf .}\label{p4.1} Assume that $c_0\leq\beta<\beta^*$, $a\geq b\beta/2$, and $(u,h)$ is a  solution of $(P)$ with $u_0=\lambda \psi$ and $\psi\in\mathscr{X}(h_0)$. Then there exist $\lambda_*\in(0,\infty]$ such that

{\rm (i)} if $0<\lambda<\lambda_*$, then vanishing happens;

{\rm (ii)} if $\lambda\geq\lambda_*$, then $h_\infty=\infty$ and $\ds\lim_{t\to\infty}u(t,x)=0$ locally uniformly in $[0,\infty)$.
\end{proposition}

Besides, we give a uniform convergence result when vanishing and virtual spreading do not happen for the solution $u$.
\begin{proposition}{\hspace{-1.8mm}\bf .}\label{p4.2} Assume that vanishing and virtual spreading do not happen for the solution $u$ of $(P)$.

{\rm (i)} if $\beta= c_0$ and $a\geq b\beta/2$, then  $\ds\lim_{t\to\infty}\|u(t,x)\|_{L^\infty([0,h(t)])}=0$;

{\rm (ii)} if $c_0<\beta<\beta^*$ and $b=0$, then $\ds\lim_{t\to\infty}\|u(t,x)-V^*(x-h(t))\|_{L^\infty([0,h(t)])}=0,$
where $V^*(z)$ is the unique solution of $(\ref{b4.9})$.
\end{proposition}
\noindent{\bf Proof.}  This proof can be accomplished successfully by following the argument of \cite[Theorem 4.15]{GLZ} if one can show
\begin{equation}\label{kkk1}\ds\lim_{t\to\infty}h'(t)=\beta-c_0\end{equation}
under the conditions for cases (i) and (ii).
The proof of this limit is divided into several steps.

{\it Step 1.} We assert that vanishing happens if $c_0<\beta<\beta^*$, $b=0$ and there exist  $t_1\geq0$ and  $x_1\in\mathbb{R}$ such that
  $$ h(t_1)<x_1,\ \,u(t_1,x)\leq V^*(x-x_1)\ \,{\rm for}\ x\in[0,h(t_1)].$$
In fact, since $V^*(x-(\beta-c_0)t-x_1)$ satisfies the first equation of $(P)$ and the free boundary condition at $x=g(t):=(\beta-c_0)(t-t_1)+x_1$, and
  $$g'(t)=\beta-c_0=-\mu(V^*)'(0),$$
we see that $(V^*(x-g(t)),g(t))$ is an upper solution of $(P)$ for $t\geq t_1$. Applying a similar argument to that of \cite[Lemma 4.10]{GLZ} gives $h_\infty<\infty$, which implies vanishing happens.

{\it Step 2.} Under the conditions for  cases (i) and (ii), if vanishing dose not happen for the solution $u$, then
\begin{equation}\label{kkkk}\lim_{t\to\infty}[h(t)-(\beta-c_0)t]=\infty.\end{equation}

If $\beta=c_0$ and $h_\infty<\infty,$ then it is easy to know that vanishing happens for $u$. This contradicts with the assumption, so (\ref{kkkk}) holds for $\beta=c_0$. For the case $c_0<\beta<\beta^*$ and $b=0$, one can use a similar discussion as in \cite[Lemma 4.12]{GLZ} to deduce that for any large $M>0$, $h(t)>(\beta-c_0)t+M$ if $t$ is sufficiently large.
Therefore (\ref{kkkk}) is obtained for the case $c_0<\beta<\beta^*$.

{\it Step 3.} Proof of (\ref{kkk1}). It suffices to show that
\begin{equation}\label{4.12}h'(t)>\beta-c_0\ \, {\rm for\  all\  large}\ t,\end{equation}  because $\ds\lim_{t\to\infty}h'(t)\leq \beta-c_0$ can be treated similarly to the proof of \cite[Lemma 4.13]{GLZ}.

It is apparent that (\ref{4.12}) is right if $\beta=c_0$. In the following we assume $c_0<\beta<\beta^*$ and $b=0$. Without loss of generality, assume
$$u_0'(0)>0, \ u_0'(h_0)<0\ {\rm and }\ u_0(x)>0\ {\rm for}\ x\in(0,h_0).$$
(Otherwise $u_0(x)$ can replaced by $u(1,x)$ in the following analysis.) So we can choose $X>h_0$ sufficiently large such that $u_0(x)$ intersects $V^*(x-M)$ at exactly two points for any $M\geq X$.   By (\ref{kkkk}), there exists $T_X>0$ such that $h(t)-(\beta-c_0)t>X$ for all $t\geq T_X$. For arbitrary $\ell>h(T_X)$, denote $T_\ell$ the unique time such that $h(T_\ell)=\ell$. Set $X_\ell:=h(T_\ell)-(\beta-c_0)T_\ell(>X)$. We can study the intersection points between $u(t,\cdot)$ and $V^*(\cdot-(\beta-c_0)t-X_\ell)$ to derive that $h'(t)>\beta-c_0$ for all $t>T_X$.
Using a similar argument as in \cite[Lemma 4.12]{GLZ},  one can show that there exists $T^*>0$ such that
 $$h'(T^*)=-\mu u_x(T^*,h(T^*))> \beta-c_0$$
and  $(\beta-c_0)t+X_\ell<h(t)\ \ {\rm for\ all}\ t>T^*$. Therefore, $T^*$ is nothing but $T_\ell$, and
$$h'(T_\ell)=-\mu u_x(T_\ell,h(T_\ell))=-\mu u_x(T_\ell,\ell)>\beta-c_0.$$
Since $\ell>h(T_X)$ is arbitrary, $T_\ell$ is continuous and strictly increasing in $\ell$, we indeed derive
$h'(t)>\beta-c_0\ \ {\rm for\ all}\ t>T_X$.
 Consequently, (\ref{4.12}) is obtained.

The proof is complete. \ \ \ \fbox{}

\noindent{\bf Proofs of Theorems \ref{t4.4} and \ref{t4.5}.} Based on Propositions \ref{p4.1} and \ref{p4.2}, we can make use of the same arguments as in \cite[Theorems 2.2 and 2.3]{GLZ} to complete the proofs of Theorems \ref{t4.4} and \ref{t4.5}.

\end{document}